\documentclass[twoside,a4paper,11pt,leqno]{amsart}

\usepackage{amssymb}
\usepackage[PostScript=dvips]{diagrams}
\swapnumbers
\makeatletter\renewcommand{\@secnumfont}{\relax}\makeatother
\theoremstyle{plain}
\newtheorem{thm}[subsection]{Theorem}
\newtheorem{prop}[subsection]{Proposition}

\newtheorem{lemma}[subsection]{Lemma}
\theoremstyle{definition}
\newtheorem{defn}[subsection]{Definition}
\newtheorem*{idefn}{Definition}
\theoremstyle{remark}
\newtheorem{rem}[subsection]{Remark}

\newcommand{\proofof}[1]{\end{#1}\begin{proof}}
\newcommand{\emphdef}{\textit}
\numberwithin{equation}{section}

\newcommand{\acknowledge}{\subsection*{Acknowledgements}}
\newcounter{numl}

\newcommand{\labelnuml}{\textup{(\roman{numl})}}
\newenvironment{numlist}{\begin{list}{\labelnuml}%
{\usecounter{numl}\setlength{\leftmargin}{0pt}%
\setlength{\itemindent}{2\parindent}%
\setlength{\itemsep}{\smallskipamount}\def
\makelabel ##1{\hss \llap {\upshape ##1}}}}{\end{list}}
\newenvironment{bulletlist}{\begin{list}{\labelitemi}%
{\setlength{\leftmargin}{\parindent}\def
\makelabel ##1{\hss \llap {\upshape ##1}}}}{\end{list}}
\newenvironment{numlproof}{\begin{proof}\begin{numlist}}
{\qed \end{numlist}\begingroup\let\qed\relax\end{proof}\endgroup}
\DeclareMathAlphabet{\mathrmsl}{OT1}{cmr}{m}{sl}
\makeatletter
\newcommand{\low}{\@ifnextchar^{}{^{\vphantom x}}}
\newcommand{\high}{\@ifnextchar_{}{_{\vphantom I}}}
\makeatother
\DeclareSymbolFont{script}{U}{eus}{m}{n}
\DeclareMathSymbol{\Wedge}{0}{script}{"5E}
\newcommand{\rssymb}[2]{\newcommand{#1}{{\mathrmsl{#2\mkern1mu}}}}
\newcommand{\calsymb}[2]{\newcommand{#1}{{\mathcal{#2}}}}
\newcommand{\bbsymb}[2]{\newcommand{#1}{{\mathbb{#2}}}}
\newcommand{\lieoper}[2]{\newcommand{#1}{\mathop{\mathfrak{#2}}}}
\newcommand{\oper}[3][n]{\newcommand{#2}{\mathop{\mathrm{#3}}\ifx
  n#1\nolimits\else\limits\fi}}
\newcommand{\rsoper}[3][n]{\newcommand{#2}{\mathop{\mathrmsl{#3\mkern1mu}}\ifx
  n#1\nolimits\else\limits\fi}}
\bbsymb\C{C}\bbsymb\HQ{H}\bbsymb\R{R}
\bbsymb\W{W}\bbsymb\V{V}\bbsymb\E{E}\bbsymb\F{F}\bbsymb\N{N}
\bbsymb\Z{Z}\bbsymb\U{U}
\calsymb\cC{C}\calsymb\cD{D}\calsymb\cE{E}\calsymb\cG{G}\calsymb\cH{H}
\calsymb\cK{K}\calsymb\cL{L}\calsymb\cN{N}\calsymb\cR{R}\calsymb\cW{W}
\newcommand{\eps}{\varepsilon}
\newcommand{\gam}{\gamma}
\newcommand{\Gam}{\Gamma}
\newcommand{\lam}{\lambda}
\newcommand{\be}{\beta}
\newcommand{\de}{\delta}
\newcommand{\si}{\sigma}
\newcommand{\al}{\alpha}
\oper\End{End}\oper\Sym{Sym}\oper\Cl{Cl}
\oper\Aut{Aut}\oper\GL{GL}\oper\SL{SL}
\oper\CO{CO}\oper\On{O}\oper\SO{SO}
\oper\CSpin{CSpin}\oper\Symp{Sp}\oper\Un{U}\oper\SU{SU}\oper\CU{CU}
\lieoper\gl{gl}\lieoper\sgl{sl}
\lieoper\co{co}\lieoper\on{o}\lieoper\so{so}
\lieoper\symp{sp}\lieoper\un{u}\lieoper\su{su}\lieoper\csu{csu}
\newcommand{\lie}[1]{\mathfrak{#1}}
\newcommand{\restr}[1]{|_{\lower.5pt\hbox{${}_{#1}$}}}

\renewcommand{\geq}{\geqslant}
\renewcommand{\leq}{\leqslant}
\newcommand{\dsum}{\oplus}               
\newcommand{\Dsum}{\bigoplus}            
\newcommand{\tens}{\mathbin{\otimes}}    
\newcommand{\subideal}{\ltimes}          
\newcommand{\intersect}{\mathinner{\cap}}
\newcommand{\act}{\mathinner{\cdot}}
\newcommand{\dual}{^{*\!}}
\newcommand{\Cinf}{\mathrm{C}^\infty}
\rsoper\kernel{ker}
\rsoper\image{im}
\rsoper\alt{alt}           
\rsoper\sym{sym}           
\rsoper\trace{tr}          
\rsoper\detm{det}          
\rsoper\divg{div}
\rsoper\Ad{Ad}
\rsoper\ad{ad}
\rssymb\iden{id}
\newcommand{\cross}{\mathbin{{\times}\!}\low}

\newcommand{\pd}{\partial}

\newcommand{\g}{\lie{g}}
\newcommand{\p}{\lie{p}}
\newcommand{\f}{\lie{f}}
\newcommand{\h}{\lie{h}}
\newcommand{\weyl}{\lie{w}}
\newcommand{\gM}{\g_M\low}
\newcommand{\pM}{\p_M\low}
\newcommand{\fM}{\f_M\low}
\newcommand{\weylM}{\weyl_M\low}
\newcommand{\T}{\lie{m}}
\newcommand{\gs}{\lie{u}}
\newcommand{\ds}{{\bullet}}
\newcommand{\cc}{\theta}                
\newcommand{\InvDer}{\nabla^\cc}        
\newcommand{\aE}{\V}\newcommand{\aF}{\W}
\newcommand{\vE}{V} \newcommand{\vF}{W} 
\newcommand{\age}{\eps\low}             
\newcommand{\Tage}{{\tilde\eps}\low}    
\newcommand{\Ge}{\cE\low}               
\newcommand{\gem}{E\low}                
\newcommand{\Tgem}{{\tilde E}\low}      
\newcommand{\sg}{\eps_\ds}                
\newcommand{\sgG}{\cE_\ds}                
\newcommand{\sgM}{E_\ds}                  

\newcommand{\MCw}{\eta}                 
\newcommand{\ve}{v}        
\newcommand{\vs}{\varphi}  
\newcommand{\vf}{f}        
\newcommand{\pt}{y}        
\newcommand{\ag}{\gam}     
\newcommand{\td}{A}        
\newcommand{\frm}{e}       
\newcommand{\dfr}{e}       
\newcommand{\rtf}{\pi_\f}  
\newcommand{\ind}{b}       
\newcommand{\x}{\boxtimes} 
\newcommand{\hj}{\hat{j}}  
\newcommand{\lroot}{\alpha}
\rsoper\project{proj}
\rsoper\gr{gr}
\newcommand\xbdbx[5]
{\begin{picture}(76,12)\put(3,3){\line(1,0){18}}%
\put(25,3){\line(1,0){6}}\put(73,3){\line(-1,0){18}}%
\put(51,3){\line(-1,0){6}}\put(39,3){\makebox(0,0){\dots}}%
\put(3,3){\makebox(0,0){$\times$}}\put(23,3){\makebox(0,0){$ \bullet$}}%
\put(53,3){\makebox(0,0){$ \bullet$}}\put(73,3){\makebox(0,0){$\times$}}%
\put(3,10){\makebox(0,0){\scriptsize $#1$}}%
\put(23,10){\makebox(0,0){\scriptsize $#2$}}%
\put(38,10){\makebox(0,0){\scriptsize $#3$}}%
\put(53,10){\makebox(0,0){\scriptsize $#4$}}%
\put(73,10){\makebox(0,0){\scriptsize $#5$}}
\end{picture}}
\begin{document}
\title[Ricci-corrected derivatives and invariant operators]
{Ricci-corrected derivatives\\ and invariant differential operators}
\author{David M. J. Calderbank}
\address{Department of Mathematics and Statistics\\ University of Edinburgh\\
King's Buildings, Mayfield Road\\ Edinburgh EH9 3JZ\\ Scotland.}
\email{davidmjc@maths.ed.ac.uk}
\author{Tammo Diemer}
\address{Mathematics Institute\\ University of Bonn\\ Beringstrasse 1\\
53113 Bonn\\ Germany.}
\curraddr{Aareal Bank, Paulinenstrasse 15,
65189 Wiesbaden. Germany.}
\email{tammo.diemer@aareal-bank.com}
\author{Vladim\smash{\'\i}r Sou\smash{\v c}ek}
\address{Mathematical Institute\\ Charles University\\ Sokolovsk\'a 83\\ Praha,
Czech Republic.}
\email{soucek@karlin.mff.cuni.cz}
\date{July 2004}
\thanks{The first author is grateful for support from the Leverhulme Trust, the
William Gordan Seggie Brown Trust and an EPSRC Advanced Research Fellowship.
The third author is supported by GA\v CR, grant Nr.\ 201/02/1390 and
MSM 113 200 007. We also thank the Erwin Schr\"odinger Institute for its
hospitality}

\begin{abstract}
We introduce the notion of Ricci-corrected differentiation in parabolic
geometry, which is a modification of covariant differentiation with better
transformation properties. This enables us to simplify the explicit formulae
for standard invariant operators given in~\cite{CSS3}, and at the same time
extend these formulae from the context of AHS structures (which include
conformal and projective structures) to the more general class of all
parabolic structures (including CR structures).
\end{abstract}
\maketitle
\section*{Introduction}

A fundamental part of differential geometry is the study of differential
invariants of geometric structures. Our concern in this paper is the explicit
construction of such invariants. More specifically, we wish to construct
invariant differential operators for a class of structures known as
\emphdef{parabolic geometries}. These geometries have attracted attention in
recent years for at least two reasons: first, they include examples of
long-standing interest in differential geometry, such as conformal structures,
projective structures and CR structures; second, they have a rich algebraic
theory, due to their intimate relation with the representation theory of
parabolic subgroups $P$ of semisimple Lie groups $G$.

A great deal of progress in our understanding of invariant differential
operators in parabolic geometry has been made through the efforts of many
people. The key idea, pioneered by Eastwood--Rice~\cite{ER} and
Baston~\cite{Baston12} is that (generalized) Bernstein--Gelfand--Gelfand (BGG)
complexes of parabolic Verma module homomorphisms are dual to complexes of
invariant linear differential operators on generalized flag varieties $G/P$,
and these complexes should admit `curved analogues', that is, there should
exist sequences of invariant linear differential operators on curved manifolds
modelled on these homogeneous spaces. This is now known to hold for all
regular BGG complexes~\cite{CSS4}.

The prototypical parabolic geometry is conformal geometry, whose model is the
generalized flag variety $\SO_0(n+1,1)/\CO(n)\subideal\R^{n*}$, namely the
$n$-sphere $S^n$ with its standard flat conformal structure.  In this case
another, more explicit, approach to the study of differential invariants has
been fruitfully applied: one first chooses a representative Riemannian metric
(or, more generally, a compatible Weyl structure); then one invokes the
well-known results of invariant theory in Riemannian (or Weyl) geometry;
finally one studies how the differential invariants depend on the choice of
Riemannian metric (or Weyl structure). The differential invariants which are
independent of the choice are invariants of the conformal structure. The
advantage of these methods is that they give explicit formulae for the
differential conformal invariants in familiar Riemannian terms.  They have
been particularly successful in the construction of first and second order
linear differential operators~\cite{Branson,Fegan,Gauduchon2}, but have been
extended to higher order operators and certain other parabolic structures
(variously known as the $|1|$-graded, abelian or AHS structures)
in~\cite{CSS3}.  In the first order case, the approach of Fegan~\cite{Fegan}
has also been extended to arbitrary parabolic geometries~\cite{SlSo}.

In this paper we build upon these explicit constructions of invariant
differential operators in terms of a compatible Weyl structure. The natural
geometries for such constructions are parabolic geometries, because we have a
good notion of Weyl structure~\cite{CSl}, similar to the conformal case, and
we can be sure than many invariant operators exist because we have the BGG
sequences~\cite{CSS4}. Our results are three-fold: first, we simplify the
formulae for standard linear differential operators given in~\cite{CSS3};
second, we extend these formulae to arbitrary parabolic structures; third, we
uncover a fundamental object in parabolic geometry, the Weyl jet operator, and
its components, the Ricci-corrected Weyl derivatives of the title of our
paper. Since it will take us a little while in the body of the text to reach
these results, we shall spend some time now explaining what the
Ricci-corrected Weyl derivatives are, in the case of conformal geometry.

Let $M$ be an $n$-dimensional manifold with a conformal structure $c$.  A
compatible Weyl connection $D$ on $M$ is a torsion-free conformal connection
on the tangent bundle of $M$. It therefore induces a connection on the
conformal frame bundle of $M$, and hence covariant derivatives on any vector
bundle $\vE$ associated to the frame bundle via a representation $\lam$ of the
conformal group $\CO(n)$ on a vector space $\aE$.

The most familiar Weyl connections are the Levi-Civita connections of
representative Riemannian metrics for $c$. However, the broader context of
Weyl connections has a few advantages in conformal geometry:
\begin{bulletlist}
\item Weyl connections form an affine space, modelled on the space of
$1$-forms $\gam$ and we write $D\mapsto D+\gam$ for this affine structure;
\item the construction of the Levi-Civita connection from a
Riemannian metric involves taking a derivative, so that differential
invariants have one order higher in the metric than in the connection;
\item a choice of Riemannian metric reduces the structure group to $\SO(n)$,
making it easy to forget the $1$-dimensional representations of the conformal
group, an omission which comes back with a vengeance in the form of conformal
weights.
\end{bulletlist}

The affine structure of the space of Weyl connections provides a
straightforward formulation of the well-known folklore that conformal
invariance only needs to be checked infinitesimally: we regard a differential
invariant $F$ constructed using a Weyl connection as a function $F(D)$; $F$ is
said to be a \emphdef{conformal invariant} if it is independent of $D$; by the
fundamental theorem of calculus, this amounts to checking that $\partial_\gam
F(D)=0$ for all Weyl connections $D$ and $1$-forms $\gam$. More generally, if
$F(D)$ is polynomial in $D$ then this dependence can be computed using
Taylor's Theorem. In particular, if $s$ is a section of an associated bundle
$\vE$ and $X$ is a vector field then $\partial_\gam D\low_Xs=[\gam,X]\act s$
where $[\gam,X]=\gam(X)\iden+\gam\wedge X\in\co(TM)=\R\,\iden_{TM}\dsum
\so(TM)$ and $\act$ denotes the natural action of $\co(TM)$ on $\vE$ (induced
by the representation $\lam$ of $\co(n)$ on $\aE$).

The Ricci-corrected derivatives have their origins in the observation that the
explicit formulae for conformally invariant differential operators in terms of
a Weyl connection appear to have a systematic form, essentially depending only
on the order of the operator. For first order operators, this is quite
straightforward~\cite{Fegan,Gauduchon2}: conformally invariant first order
linear operators are all of the form $\pi\circ
D\colon\Cinf(M,\vE)\to\Cinf(M,\vF)$, where $\vE$ and $\vF$ are associated
bundles, $D$ the covariant derivative on $\vE$ induced by a Weyl connection,
and $\pi$ is induced by an equivariant map $\R^{n*}\tens\aE\to\aF$.  Evidently
$\partial_\gam \pi(Ds)=\pi\bigl([\gam,\cdot]\act s\bigr)$, so we obtain
conformally invariant operators by letting $\pi$ be the projection onto the
zero eigenspace of the operator $\Psi\in\End( T\dual M\tens\vE)$ defined by
$\Psi(\gam\tens s)=[\gam,\cdot]\act s$. Since these eigenvalues can be
shifted, by tensoring $\V$ with a one dimensional representation of $\CO(n)$,
a large number of first order operators are obtained.

Ricci corrections make their first appearance at the level of second order
operators. Here one finds that many conformally invariant operators are of the
form $\pi(D^2s+r^D\tens s)$ where $\pi$ is induced by an equivariant map
$\R^{n*}\tens\R^{n*}\tens\aE\to S^2\R^{n*}\tens\aE\to \aF$ (the first map
being symmetrization), and $r^D$ is the \emphdef{normalized Ricci curvature}
of the Weyl connection, which is a covector-valued $1$-form on $M$ constructed
from the curvature of the Weyl connection. For present purposes, all we need
to know about $r^D$ is its dependence on $D$: $\partial_\gam r^D=-D\gam$.

Now compare this with the variation of the second derivative:
\begin{align*}
\partial_\gam D^2_{X,Y}s&=[\gam,X]\act D\low_Ys-D\low_{[\gam,X]\act Y}s
+D\low_X\bigl([\gam,Y]\act s\bigr)-[\gam,D\low_XY]\act s\\
&=[D\low_X\gam,Y]\act s+[\gam,X]\act D\low_Ys
+[\gam,Y]\act D\low_Xs-D\low_{[\gam,X]\act Y}s.
\end{align*}
This formula means that we can use the Ricci curvature to make the second
derivative \emph{algebraic} in $D$.

\begin{idefn} The \emphdef{Ricci-corrected second derivative} on sections $s$
of an associated bundle $\vE$ is defined by
$D^{(2)}_{X,Y}s=D^2_{X,Y}s+[r^D(X),Y]\act s$.
\end{idefn}
Hence $\partial_\gam D^{(2)}s=\gam\ast_1 Ds$ where
\begin{align*}
(\gam\ast_1\phi)\low_{X,Y}&=
[\gam,X]\act\phi\low_Y+[\gam,Y]\act\phi\low_X-\phi\low_{[\gam,X]\act Y}\\
&=([\gam,X]\act\phi)\low_Y+[\gam,Y]\act\phi\low_X.
\end{align*}

It is now a purely algebraic matter to find projections $\pi$ such that $\pi
\bigl(D^{(2)}s\bigr)$ is a conformal invariant of $s$. As it turns out, these
projections often have the property that $\pi\bigl([r^D(\cdot),\cdot]\act
s\bigr)=\pi(r^D\tens s)$.  The simplest example is the conformal hessian
$\sym_0D^{(2)}s=\sym_0 (D^2s+r^Ds)$ where $s$ is a section of the weight $1$
line bundle $L$.

The same ideas apply to higher order operators: we want to write these
operators as $\pi\circ D^{(k)}$ for some projections induced by an equivariant
map $(\tens^k \R^{n*})\tens\V\to S^k\R^{n*}\tens\V\to\W$, where $D^{(k)}$ is a
Ricci-corrected $k$th power of the Weyl connection.  Again we make the
observation that if $\partial_\gam(\pi\circ D^{(k)})=0$ then $\pi\circ
D^{(k)}$ must certainly be algebraic in $D$, and we can in fact arrange for
$D^{(k)}$ itself to be algebraic in $D$.

\begin{idefn} The \emphdef{Ricci-corrected powers of the Weyl connection} on
an associated bundle $\vE$ are defined inductively by $D^{(0)}s=s$,
$D^{(1)}s=Ds$ and
\begin{equation*}
\iota\low_X D^{(k+1)}s=D\low_XD^{(k)}s+r^D(X)\ast_1 D^{(k-1)}s.
\end{equation*}
where
\begin{equation*}
(\gam\ast_1\phi)_{X_1,\ldots X_k}
=\sum_{j=1}^k([\gam,X_i]\act\phi\low_{X_1,\ldots
X_{i-1}})\low_{X_{i+1},\ldots X_k}.
\end{equation*}
A calculation shows that $\pd_\gam D^{(k)}s=\gam\ast_1 D^{(k-1)}s$.
\end{idefn}
The inductive formula can easily be summed to give
\begin{multline*}\notag
D^{(k)}= \sum_{\ell+m=k}\;
\sum_{1\leq i_1<\cdots<i_\ell\leq m}
D^{m-i_\ell}\circ (r^D(\cdot)\ast_1)\circ D^{i_\ell-i_{\ell-1}-1}\circ\cdots\\
\cdots\circ (r^D(\cdot)\ast_1)\circ D^{i_2-i_1-1}\circ(r^D(\cdot)\ast_1)
\circ D^{i_1-1}.
\end{multline*}
Thus the search for explicit invariant operators reduces to an algebraic
problem, and we shall find that a large class of projections $\pi$
annihilating $\gam\ast_1 D^{(k-1)}s$ produce universal numbers when applied
to the terms of $D^{(k)}$, essentially because the action $\gam\ast_1$ on
$(\tens^j T\dual M)\tens\vE$ for $j<k$ is closely related to the action on
$(\tens^k T\dual M)\tens\vE$.

This theory of Ricci-corrected derivatives in conformal geometry was developed
over several years by the first two authors, and described, in part,
in~\cite{Thesis}. In the homogeneous case, i.e., on $S^n$, the second author
explained the action $\gam\ast_1$ in terms of the second order part of the
action of a conformal vector field on sections of a homogeneous vector
bundle~\cite{Thesis}.  It became clear however, that there was a systematic
underlying principle behind these formulae, even in the curved case, which
should also generalize to arbitrary parabolic geometries. More precisely, the
action $\gam\ast_1$ is related to part of the action of the nilradical of the
parabolic subalgebra $\p$ on certain \emphdef{semiholonomic jet modules} $\hat
J^k_0\aE$ associated to a module $\aE$~\cite{ES,Slovak}. This action is
considerably more complicated in general than it is in the conformal case, and
we are forced to consider the representation theory of the entire parabolic
subgroup $P$, not just its Levi factor $P_0$ as we did in the conformal case
(where $P_0=\CO(n)$ and $P=\CO(n)\subideal\R^{n*}$). Hence we must work with
$P$-modules $\aE$, and the corresponding vector bundles $\vE$ are associated
to a larger principal $P$-bundle, which in the conformal case is the Cartan
bundle with its normal Cartan connection. This development ultimately provides
a simple conceptual explanation for the formulae we obtain.

The structure of the paper is as follows. We begin by defining Cartan
geometries and invariant differentiation in section~\ref{idcg}: this is a
standard way to treat parabolic geometries~\cite{CD,CG,Slovak}, although in
practice a geometry is defined by more primitive data, which must be
differentiated to obtain the Cartan connection~\cite{CS}. The key point is
that a Cartan connection determines semiholonomic jet operators taking values
in bundles associated to semiholonomic jet modules.

In sections~\ref{ps}--\ref{pgws} we begin the study of parabolic geometries
and Weyl structures. We adopt a novel approach to Weyl structures (which we
relate to the approach of \v Cap and Slovak~\cite{CSl} in Appendix A) in order
to emphasise the relationship between the geometry of Weyl structures and some
elementary representation theory which we exploit throughout our treatment. At
the algebraic level, a Weyl structure is a lift $\age$ of a certain `grading
element' $\age_0$ in the Levi factor $\p_0$ to the parabolic Lie algebra
$\p$. The grading element induces a filtration of a $P$-module $\aE$ and a
lift $\age$ splits this filtration, i.e., determines an isomorphism $\age_\aE$
of $\aE$ with its associated graded module $\gr\aE$. Geometric Weyl structures
for parabolic geometries are given simply by applying the same procedure
pointwise on the underlying manifold: a geometric Weyl structure $\gem$ then
determines a splitting $\gem_\vE\colon \vE\to \gr\vE$ of any associated
filtered $P$-bundle $\vE$.  Now if $\vE$ is a filtered $P$-bundle, so is its
semiholonomic $k$-jet bundle $\hat J^k\vE$, and splitting this bundle allows
us to project out Ricci corrected derivatives as components of the
$k$-jet. This simple construction, which we present in section~\ref{rcd},
makes Ricci corrected differentiation easy to study from a theoretical point
of view.  On the other hand, it is also easily related to covariant
differentiation: explicit formulae are obtained as soon as one understands the
action of the nilradical of $\p$ on jet modules.

In section~\ref{sio}, we pave the way for the construction of invariant
operators by studying special types of projections from jet modules, which
have the effect of killing most of the complicated terms in the jet module
action. For irreducible modules, the remaining terms of the jet module action
reduce to the projection of a scalar action, which we compute using
Casimirs. In section~\ref{ecio} we give our construction of a large class of
invariant operators, and write out the formulae for operators up to order
$8$. We illustrate the scope of the constructions in section~\ref{scope}
and give some examples in conformal geometry.

Finally let us mention further potential applications of Ricci-corrected
differentiation.  Although in this paper we have applied Ricci-corrected
derivatives to the construction of invariant \emphdef{linear} differential
operators, the same ideas can be expected to yield explicit formulae for
multilinear differential operators, such as the operators
of~\cite{CD}. Indeed, this was our original motivation to study
Ricci-corrected differentiation in conformal geometry: one approach to
construct (say) bilinear differential operators is to combine terms
constructed from pairs of \emph{noninvariant} linear differential operators;
to do this one needs noninvariant operators which nevertheless depend on the
choice of Weyl structure  in a simple way---projections of Ricci-corrected
derivatives onto irreducible components have this property.

\acknowledge We would like to thank Fran Burstall and Jan Slovak for helpful
discussions concerning Weyl structures.

\section{Invariant derivatives in Cartan geometry}\label{idcg}

Parabolic geometries are geometries modelled on a generalized flag variety
$G/P$, i.e., $G$ is a semisimple Lie group and $P$ is a parabolic subgroup.  A
standard way to define `curved versions' of homogeneous spaces is as Cartan
geometries. In this section, we recall the basic calculus of such geometries,
following~\cite{CD,CG,ES,Sharpe,Slovak}.

Fix a Lie algebra $\g$ with a Lie group $P$ acting by automorphisms such that
$\p$ is a $P$-equivariant subalgebra of $\g$, and the derivative of the
$P$-action on $\g$ is the adjoint action of $\p$ on $\g$.  (These technical
conditions are simply those that arise when $P$ is a subgroup of a Lie group
$G$ with Lie algebra $\g$.)

\begin{defn}
Let $M$ be a manifold with dimension $\dim M=\dim \g-\dim \p$.  A
\emphdef{Cartan connection} of type $(\g,P)$ on $M$ is a principal $P$-bundle
$\pi\colon \cG\to M$, together with a $P$-invariant $\g$-valued $1$-form
$\cc\colon T\cG\to\g$ such that for each $\pt\in\cG$, $\cc_\pt\colon
T_\pt\cG\to\g$ is an isomorphism which sends each generator $\zeta_\xi$ of the
$P$-action to the corresponding $\xi\in\p$, i.e.,
$\zeta_{\xi,\pt}=\cc^{-1}_\pt(\xi)$. Here $P$-invariance means that
$\Ad(p)\act r_p^*\cc=\cc$ for any $p\in P$, where $r_p$ denotes the right
$P$-action on $\cG$. We refer to $(M,\cG,\cc)$ as a \emphdef{Cartan geometry}.
\end{defn}
Note that Cartan connections form an \emph{open subset} of an affine space
modelled on the space of horizontal $P$-invariant $\g$-valued $1$-forms.
However, one can freely add horizontal $P$-invariant \emph{$\p$-valued}
$1$-forms to a Cartan connection, without losing invertibility.

Associated to any $P$-module $\aE$ is a vector bundle $\vE=\cG\cross_P\aE$,
defined to be the quotient of $\cG\cross\aE$ by the action $(\pt,\ve)\mapsto
(\pt p^{-1},p\act\ve)$. This induces an action $(p\act \vf)(\pt)=p\act \vf(\pt
p)$ on functions $\vf\in\Cinf(\cG,\aE)$ which identifies sections $\vs$ of
$\vE=\cG\cross_P\aE$ over $M$ with $P$-invariant functions
$\vf\colon\cG\to \aE$:
\begin{equation*}\notag
\Cinf(M,\vE)=\Cinf(\cG,\aE)^P.
\end{equation*}
Similarly $P$-invariant horizontal $\aE$-valued forms on $\cG$ are identified
with forms on $M$ with values in $\vE$. In particular the $1$-form $T\cG\to
\g/\p$ induced by the Cartan connection $\cc$ is $P$-invariant and horizontal,
corresponding to a bundle map $TM\to \cG\cross_P\g/\p$.  The open condition on
the Cartan connection means that this is an isomorphism and henceforth we
identify $TM$ with $\cG\cross_P\g/\p$ in this way. We let $\gM=\cG\cross_P\g$
and observe that there is a surjective bundle map from $\gM$ to $TM$, with
kernel $\pM=\cG\cross_P\p$.

Cartan connections do not in general induce covariant derivatives on
associated bundles, but there is a way of differentiating sections of such
bundles using $\gM$ instead of $TM$.

\begin{defn} Let $(\cG,\cc)$ be a Cartan connection of type $(\g,P)$ on
$M$, and let $\aE$ be a $P$-module with associated vector bundle
$\vE=\cG\cross_P \aE$.  Then the linear map defined by
\begin{align*}
\InvDer\colon\Cinf(\cG,\aE)&\to\Cinf(\cG,\g\dual\tens\aE)\\
\InvDer_\xi \vf &= d\vf\bigl(\cc^{-1}(\xi)\bigr)
\end{align*}
(for all $\xi$ in $\g$) is $P$-equivariant. The restriction to
$\Cinf(\cG,\aE)^P$, or equivalently the induced linear map $\InvDer\colon
\Cinf(M,\vE)\to\Cinf(M,\gM^*\tens \vE)$, is called the \emphdef{invariant
derivative} on $\vE$.
\end{defn}

The \emphdef{curvature} $K\colon\Wedge^2T\cG\to\g$ of a Cartan geometry is
defined by
\begin{equation*}
K(X,Y)=d\cc(X,Y)+[\cc(X),\cc(Y)].
\end{equation*}
It induces a curvature function
$\kappa\colon\cG\to\Wedge^2\g\dual\tens\g$ via
\begin{equation*}
\kappa_\pt(\xi,\chi)=K_\pt\bigl(\cc^{-1}(\xi),\cc^{-1}(\chi)\bigr)=
[\xi,\chi]-\cc_\pt[\cc^{-1}(\xi),\cc^{-1}(\chi)],
\end{equation*}
where $\pt\in\cG$ and the latter bracket is the Lie bracket of vector fields
on $\cG$.

The $\p$-invariance of $\chi\mapsto\cc^{-1}(\chi)$ for $\chi\in\g$ means that
$[\cc^{-1}(\xi),\cc^{-1}(\chi)]=\cc^{-1}[\xi,\chi]$ for any $\xi\in\p$, and
hence that $\kappa(\xi,\cdot)=0$ for $\xi\in\p$ so that
$\kappa\in\Cinf(\cG,\Wedge^2(\g/\p)\dual\tens\g)^P$.  In other words $K$ is a
horizontal $2$-form and induces $K\low_M\in\Cinf(M,\Wedge^2T\dual M\tens\gM)$.

\begin{lemma} Let $(\cG,\cc)$ be a Cartan geometry of type $(\g,P)$ on $M$.
\begin{numlist}
\item For $\vf\in\Cinf(\cG,\aE)^P$, we have $(\InvDer_\xi \vf)(\pt)
+\xi\act\bigl(\vf(\pt)\bigr)=0$ for all $\xi\in\p$ and $\pt\in\cG$.
\item We also have $\InvDer_\xi(\InvDer_\chi \vf)-\InvDer_\chi(\InvDer_\xi
\vf) = \InvDer_{[\xi,\chi]}\vf - \InvDer_{\kappa(\xi,\chi)}\vf$ for all
$\xi,\chi\in\p$.
\end{numlist}
\end{lemma}
\begin{numlproof}
\item Differentiate the $P$-invariance condition $p\act(\vf(\pt p))=\vf(\pt)$.
\item Both sides are equal to $d\vf([\cc^{-1}(\xi),\cc^{-1}(\chi)])$.
\end{numlproof}

These facts enable us to define a semiholonomic jet operator $\hj^k_\cc$
identifying the semiholonomic jet bundle $\hat J^k\vE$ with an associated
bundle $\cG\cross_P\hat J^k_0\aE$~\cite{CD,CSS1,Slovak}.  Recall that the
semiholonomic jet bundles are defined inductively by $\hat J^1\vE=J^1\vE$ and
$\hat J^{k+1}\vE$ is the subbundle of $J^1\hat J^k\vE$ on which the two
natural maps to $J^1\hat J^{k-1}\vE$ agree. The advantage of semiholonomic
jets is that they depend only on the $1$-jet functor, the natural
transformation $J^1\vE\to \vE$ and some abstract nonsense.

\begin{prop} Let  $(\cG,\cc)$ be a Cartan geometry of type $(\g,P)$ on $M$
and $\aE$ a $P$-module.
\begin{numlist}
\item The map $j^1_\cc\colon\Cinf(M,\vE)\to\Cinf\bigl(M,\vE\dsum(\gM^*\tens
\vE)\bigr)$ sending $\vs$ to $(\vs,\InvDer\vs)$ defines an injective bundle
map, from the $1$-jet bundle $J^1\vE$ to $\vE\dsum \bigl(\gM^*\tens\vE\bigr)$,
whose image is $\cG\cross_P J^1_0\aE$ where
$J^1_0\aE=\{(\phi_0,\phi_1)\in\aE\dsum\bigl(\g\dual\tens
\aE\bigr):\phi_1(\xi)+\xi\act\phi_0=0 \textup{ for all } \xi\in\p\}$.
\item Similarly the map $\hj^k_\cc$ sending a section $\vs$ to
$\bigl(\vs,\InvDer \vs,(\InvDer)^2\vs,\ldots(\InvDer)^k\vs\bigr)$
defines an isomorphism between the semiholonomic jet bundle $\hat J^k\vE$
and the subbundle $\cG\cross_P \hat J^k_0\aE$ of
$\Dsum_{j=0}^k\bigl((\tens^j\gM^*)\tens\vE\bigr)$,
where $\hat J^k_0\aE$ is the set of all $(\phi_0,\phi_1,\ldots \phi_k)$ in
$\Dsum_{j=0}^k\bigl((\tens^j\g\dual)\tens\aE\bigr)$ satisfying \textup(for
$1\leq i<j\leq k$\textup) the equations
\begin{align*}
\quad\phi_{j}(\xi_1,\ldots\xi_i,\xi_{i+1},\ldots\xi_j)
-\phi_{j}(\xi_1,\ldots\xi_{i+1},\xi_i,\ldots\xi_j)
&=\phi_{j-1}(\xi_1,\ldots[\xi_i,\xi_{i+1}],\ldots\xi_j)\\
\phi_{i}(\xi_1,\ldots\xi_i)
+\xi_i\act\bigl(\phi_{i-1}(\xi_1,\ldots\xi_{i-1})\bigr)&=0
\end{align*}
for all $\xi_1,\ldots\xi_j\in\g$ with $\xi_i\in\p$.
\end{numlist}
\end{prop}
\begin{numlproof}
\item Certainly the map on smooth sections only depends on the $1$-jet at each
point, and it is injective since the symbol of $\InvDer$ is the inclusion
$T\dual M\tens\vE\to \gM^*\tens\vE$. It maps into $\cG\cross_P J^1_0\aE$ by
vertical triviality, but this has the same rank as $J^1\vE$.
\item Similar: the equations are those given by the vertical triviality and
the Ricci identity, bearing in mind that $\kappa$ is horizontal. The
(semiholonomic) symbols of the iterated invariant derivatives are still
given by inclusions $(\tens^j T\dual M)\tens\vE\to(\tens^j\gM^*)\tens
\vE$.
\end{numlproof}

\section{Parabolic subalgebras and algebraic Weyl structures}\label{ps}

Parabolic geometries are Cartan geometries of type $(\g,P)$ where $\g$ is
semisimple and the Lie algebra $\p$ of $P$ is a parabolic subalgebra. In this
section we develop a few basic facts about parabolic subalgebras and their
representations, emphasising the relation between filtered and graded modules.
The key feature of parabolic subalgebras is the presence of `algebraic Weyl
structures' which split filtered modules. Such splittings, carried out
pointwise, equip parabolic geometries with covariant derivatives on associated
bundles.

A parabolic subalgebra of a semisimple Lie algebra is a subalgebra containing
a Borel (i.e., maximal solvable) subalgebra. However, to keep our treatment as
self-contained as possible, with minimal use of structure theory, we find the
following equivalent and elementary definition more convenient. We refer
to~\cite{BE,CSl,Slovak} for an alternative approach.

\begin{defn} Let $\g$ be a semisimple Lie algebra. For a subspace $\gs$
of $\g$ we let $\gs^\perp$ be the orthogonal subspace with respect to the
Killing form $(\cdot\,,\cdot)$. Then a subalgebra $\p$ of $\g$ is
\emphdef{parabolic} iff $\p^\perp$ is the \emphdef{nilradical} of $\p$,
i.e., its maximal nilpotent ideal.  It follows that the quotient
$\p_0:=\p/\p^\perp$ is a reductive Lie algebra, called the \emphdef{Levi
factor}.
\end{defn}

Let $\p^\perp, (\p^\perp)^2=[\p^\perp,\p^\perp], \ldots (\p^\perp)^{j+1}
=[\p^\perp,(\p^\perp)^j],\ldots$ be the descending central series of
$\p^\perp$.  Since $\p^\perp$ is nilpotent there is an integer $k\geq 0$,
called the \emphdef{depth} of $\p$, such that $(\p^\perp)^{k+1}=0$ but
$(\p^\perp)^k\neq 0$. Thus $\p^\perp$ has a $k$-step filtration ($k=0$ is the
trivial case $\p^\perp=0$ and $\p=\g$). We obtain from this a filtration of
$\g$ by setting $\g_{(-j)}=(\p^\perp)^j$ and $\g_{(j-1)}=\g_{(-j)}^\perp$ for
$j\geq 1$ so that
\begin{equation*}\notag
\g=\g_{(k)} \supset \g_{(k-1)} \supset \cdots \supset \g_{(1)} \supset
\g_{(0)}=\p \supset \p^\perp=\g_{(-1)} \supset \cdots \supset \g_{(-k)}
\supset \g_{(-k-1)} = 0.
\end{equation*}
It is easily verified that $[\g_{(i)},\g_{(j)}]\subseteq \g_{(i+j)}$, so that
$\g$ is a filtered Lie algebra. The associated graded Lie algebra is
$\gr\g=\Dsum_{j=-k}^k \g_j$ where $\g_j=\g_{(j)}/\g_{(j-1)}$ and is said to be
$|k|$-graded. Note in particular that $\p_0=\g_0$.

An important fact about parabolic subalgebras is that this filtration of $\g$
is split.
\begin{lemma} There are \textup(non-canonical\textup) splittings of the exact
sequences
\begin{equation}\label{es}
0\to \g_{(j-1)} \to \g_{(j)}\to \g_j\to 0
\end{equation}
which induce a Lie algebra isomorphism between $\g$ and $\gr \g$. 
\end{lemma}
\begin{proof} Any semisimple Lie algebra admits a Cartan involution, i.e.,
an automorphism $\sigma\colon\g\to\g$ such that $\sigma^2=\iden$ and
$h(\xi,\chi):=(\sigma(\xi),\chi)$ is positive definite. We split~\eqref{es}
for each $j$ by identifying $\g_j$ with the $h$-orthogonal complement to
$\g_{(j-1)}$ in $\g_{(j)}$. Suppose $\xi\in\g_{(i)}$ is $h$-orthogonal to
$\g_{(i-1)}$, i.e., $\sigma(\xi)\in\g_{(i-1)}^\perp=\g_{(-i)}$, and
$\chi\in\g_{(j)}$ is $h$-orthogonal to $\g_{(j-1)}$. Then $[\xi,\chi]\in
\g_{(i+j)}$ and $\sigma[\xi,\chi]= [\sigma(\xi),\sigma(\chi)]\in \g_{(-i-j)}$
so $[\xi,\chi]$ is $h$-orthogonal to $\g_{(i+j-1)}$. Hence the splittings
defined by $\sigma$ induce a Lie algebra isomorphism.
\end{proof}

We refer to such a splitting of $\g$ as an \emphdef{algebraic Weyl structure}:
it is not unique, but we can obtain very good control over the possible
splittings thanks to the following.
\begin{lemma}
There is a unique element $\age_0$ in the centre of $\p_0=\p/\p^\perp$ such
that $[\age_0,\xi]=j \xi$ for all $\xi\in \g_j$ and all $j$.
\end{lemma}
\begin{proof}
Since $\gr\g$ is semisimple, the derivation defined by $\xi\mapsto j \xi$ for
$\xi\in\g_j$ must be inner, i.e., equal $\ad\age_0$ for $\age_0\in\gr\g$
(which is unique since $Z(\g)=0$). Now $[\age_0,\age_0]=0$ and
$[\age_0,\xi]=0$ for all $\xi\in \g_0$, so $\age_0$ is in the centre of
$\g_0=\p_0$.
\end{proof}

\begin{defn} The element $\age_0$ is called the \emphdef{grading element}.
Let $\weyl=\{\age\in \p:\pi_0(\age)=\age_0\}$ be the set of all lifts of
$\age_0$ to $\p$ with respect to the exact sequence
\begin{equation}
0\to \p^\perp \to \p \xrightarrow{\pi_0} \p_0 \to 0.
\end{equation}
\end{defn}
The elements of $\weyl$ are precisely the algebraic Weyl structures: the
isomorphism of $\g$ with $\gr\g$ is given by the eigenspace decomposition of
$\ad\age$ for a lift of $\age_0$ to $\age\in\p\subseteq\g$.  The space of
algebraic Weyl structures is therefore $\weyl$, an affine space modelled on
$\p^\perp$.

Let $P$ be a Lie group acting on $\g$ with Lie algebra $\p$ as in the previous
section, and suppose additionally that the quotient group $P_0=P/\exp\p^\perp$
stabilizes $\age_0\in\p_0$ (which is automatic if $P_0$ is connected) so that
the adjoint action of $P$ on $\p$ preserves $\weyl$.
\begin{lemma} $\exp\p^\perp\leq P$ acts freely and transitively on $\weyl$.
\end{lemma}
\begin{proof}
If $\ag\in\p^\perp$, $(\Ad\exp\ag)\age=\exp(\ad\ag)\age
=\age+[\ag,\age]+\cdots$. The result follows because $\ad\ag$ is nilpotent on
$\p$, and $\ag\mapsto[\ag,\age]$ is a bijection on $\p^\perp$.
\end{proof}
The stabilizer of $\age$ is thus a subgroup of $P$ projecting isomorphically
onto $P_0$, so that an algebraic Weyl structure splits the quotient group
homomorphism $\pi_0\colon P\to P_0$.

The fundamental vector fields $\zeta_\ag$ ($\ag\in\p^\perp$) generating the
action of $\exp\p^\perp$ on $\weyl$ give rise to a Maurer--Cartan form
$\MCw\colon T\weyl\to \p^\perp$ with $\MCw(\zeta_\ag)=\ag$.  If we identify
$T\weyl$ with $\weyl\times\p^\perp$ using the affine space structure then
$\zeta_{\ag,\age}=[\ag,\age]$, so $\MCw_{\age}$ is the inverse of
$\ag\mapsto[\ag,\age]$ on $\p^\perp$.

\subsection{Filtered and graded modules}

We say that a $P_0$-module is \emphdef{semisimple} if it is completely
reducible and the grading element $\age_0$ acts by a scalar on irreducible
components (the latter condition is automatic for complex modules). The
eigenvalues of $\age_0$ will be called the \emphdef{geometric weights} of the
module.

There is a one to one correspondence between $P_0$-modules and $P$-modules on
which $\exp\p^\perp$ acts trivially. We say such a $P$-module is semisimple if
the corresponding $P_0$-module is.  More generally, we shall consider
$P$-modules $\aE$ with a $P$-invariant filtration
\begin{equation}\label{repfilt}
\aE=\aE_{(\lambda)}\supset
\aE_{(\lambda-1)}\supset\aE_{(\lambda-2)}\supset\cdots \supset
\aE_{(\lambda-\ell)} \supset 0
\end{equation}
(for a scalar $\lambda$ and an integer $\ell$) such that the associated graded
module $\gr\aE$ is a semisimple $P$-module graded by geometric weight. We
refer to such modules as \emphdef{filtered $P$-modules}. We extend the
definition in a straightforward way to direct sums of such modules.

An algebraic Weyl structure $\age$ splits any filtered $P$-module $\aE$ into
the eigenspaces of $\age$, giving a vector space isomorphism $\age_\aE\colon
\aE\to \gr\aE$. This isomorphism is not $P$-equivariant (though it is
tautologically $P_0$-equivariant using the splitting of $\pi_0\colon P\to P_0$
defined by $\age$). However, by the naturality of the construction, the map
$\age\mapsto \age_\aE$ is $P$-equivariant. More precisely, if $\Tage = (\Ad
p)\age$ is any other algebraic Weyl structure (with $p\in P$), then
\begin{equation}\label{gr-change}
\Tage_\aE(v)=p\act \age_\aE(p^{-1}\act v).
\end{equation}
In particular, if $p\in\exp\p^\perp$, we have $\Tage_\aE = \age_\aE\circ
(q^{-1}\act{})$.  Since the action of $\exp\p^{\perp}$ on $\weyl$ is free and
transitive, this dependence implies that if we have a smooth family
$\age=\age(s)$ of algebraic Weyl structures then
$d\age_\aE(X)(v)=-\age_\aE(\MCw(d\age(X))\act v)$.

\begin{lemma} \label{MC} For a smooth map $\eps\colon S\to \weyl$
and a filtered $P$-module $\aE$, the $\End\aE$-valued $1$-form
$\age^{-1}_\aE d\age_\aE$ on $S$ is given by the natural action of
$-\age^*\MCw$ on $\aE$.
\end{lemma}

\subsection{The tangent module}

We end this section by considering the module $\T:=\g/\p$, which is the
filtered $P$-module dual to $\p^\perp$. More precisely, the Killing form of
$\g$ gives a nondegenerate pairing between $\g/\p$ and $\p^\perp$, which will
also be denoted by $\T\dual$.  This duality depends on the normalization of
the Killing form, which we do not specify at present.

We have seen that $\g$ is a filtered $P$-module, with $\T$ and $\T\dual$ as
quotient and sub- modules respectively.  The associated graded modules $\gr\T$
and $\gr\T\dual$ are graded nilpotent subalgebras of $\gr\g$.  In particular,
as semisimple $P$-modules, $\gr\g = \gr\T\dsum\p_0\dsum \gr\T\dual$, although
the Lie bracket is not compatible with this direct sum decomposition. An
algebraic Weyl structure $\age$ therefore determines a vector space
isomorphism $\sg\colon \T\dsum\p_0\dsum \T\dual \to \g$.

Observe that $\f:=\g_1$ is the lowest geometric weight subspace of $\gr\T$,
and so is a $P$-submodule of $\T$; the dual $\f\dual$ is naturally a quotient
$P$-module of $\p^\perp=\T\dual$.

\section{Parabolic geometries and Weyl structures}\label{pgws}

\begin{defn} A \emphdef{parabolic geometry} on $M$ is a Cartan geometry
$(\cG,\cc)$ of type $(\g,P)$ with $\g$ semisimple and $\p$ parabolic,
satisfying the conditions of the previous two sections.

We define $\cG_0$ to be the principal $P_0$-bundle $\cG/\exp\p^\perp$ and we
let $\pi_0$ also denote the projection $\cG\to\cG_0$, so that $\pi_0(\pt p)
=\pi_0(\pt)\pi_0(p)$ for $\pt\in\cG$ and $p\in P$.
\end{defn}

The tangent bundle $TM=\cG\cross_P\T$ has a natural filtration induced by the
filtration of $\T$, the smallest nontrivial distribution in the filtration
being $\fM=\cG\cross_P\f$. The cotangent bundle $T\dual
M=\cG\cross_P\T\dual=\cG\cross_P\p^\perp$ is a bundle of nilpotent Lie
algebras, the nilradical bundle of $\pM=\cG\cross_P\p$.  The quotient
$\pM/T\dual M$ is a reductive Lie algebra bundle, namely
$\p\low_{M,0}:=\cG\cross_P\p_0$. Observe that $\p\low_{M,0}$ has a canonical
grading section $\gem_0$, induced by the grading element $\age_0$ of $\p_0$,
which is $P$-invariant.

\begin{defn}
Let $(\cG,\cc)$ be a parabolic geometry on $M$. Then a (geometric)
\emphdef{Weyl structure} $\gem$ on $M$ is a smooth lift of the grading section
$\gem_0$ to a section of $\pM$.
\end{defn}
Thus a Weyl structure amounts to a smooth choice of algebraic Weyl structure
at each point. Since algebraic Weyl structures form an affine space, a Weyl
structure is a section of an affine bundle, the \emphdef{bundle of Weyl
geometries} $\weylM=\cG\cross_P \weyl$. In particular, Weyl structures always
exist, and form an affine space modelled on the space of $1$-forms on $M$.

\begin{rem}[Key observation]\label{keyrem}
Any construction with algebraic Weyl structures can be carried out with
geometric Weyl structures. We can either work with associated bundles or on
the principal bundle $\cG$, and both points of view are useful.
\begin{numlist}
\item \label{kr1} If $\vE=\cG\cross_P \aE$ is bundle associated to a filtered
$P$-module $\aE$, with graded bundle $\gr\vE=\cG\cross_P\gr\aE=
\cG_0\cross_{P_0}\gr\aE$, then a Weyl structure $\gem$ provides an isomorphism
$\gem_\vE\colon \vE\to \gr\vE$, simply by applying the construction of the
previous section pointwise. We also obtain a bundle isomorphism $\sgM\colon
TM\dsum \p\low_{M,0}\dsum T\dual M\to \gM$.
\item \label{kr2} A Weyl structure may equally be regarded as a $P$-invariant
function $\Ge\colon \cG\to \weyl$. For any filtered $P$-module $\aE$, we then
have a $P$-equivariant isomorphism $\Ge_{\aE}\colon \cG\cross\aE\to \cG\cross
\gr\aE$ whose fibre at $\pt\in\cG$ is $\Ge(\pt)\low_\aE$; the induced
isomorphism of associated bundles is $\gem_\vE$.  Similarly we get a
$P$-equivariant isomorphism $\sgG\colon \cG\cross(\T\dsum \p_0\dsum
\T\dual)\to\cG\cross\g$ inducing $\sgM$.
\end{numlist}
\end{rem}

If we fix an algebraic Weyl structure then any geometric Weyl structure may be
written $\Ge=(\Ad q)\age$, where $q\colon\cG\to\exp\p^\perp$ is $P$-invariant
in the sense that $pq(\pt p)\pi_0(p^{-1})=q(\pt)$: here $P_0$ acts on $P$ via
the lift defined by $\age$.  As we discuss in Appendix A, this allows us to
relate our approach to Weyl structures to the original approach of \v
Cap--Slov\'ak~\cite{CSl}.

\section{Ricci-corrected Weyl differentiation}\label{rcd}

The main difficulty in the study of invariant differential operators on
parabolic geometries is that there is no natural covariant derivative on
associated bundles: we only have the invariant derivative $\InvDer
\colon\Cinf(M,\vE)\to\Cinf(M,\gM^*\tens \vE)$ in general. There is no
canonical projection $\gM^*\to T\dual M$; equivalently, the restriction map
$\g\dual\to\T\dual=\p^\perp$ is not $P$-equivariant.

Weyl structures provide two solutions to this problem, one well known (Weyl
connections), the other implicitly known, but not properly formalized
(Ricci-corrected Weyl connections). In our theory, both can be defined
straightforwardly using Remark~\ref{keyrem}\ref{kr1}.

\begin{defn} Let $(\cG\to M,\cc)$ be a parabolic geometry and $\gem$ be
a Weyl structure on $M$. Let $\vE=\cG\cross_P\aE$ be a filtered $P$-bundle
(i.e., associated to a filtered $P$-module).
\begin{numlist}
\item The \emphdef{Ricci-corrected Weyl connection} $D^{(1)}\colon
\Cinf(M,\vE) \to \Cinf(M, T\dual M\otimes \vE)$ is given by $D^{(1)}_X \vs=
\InvDer_{\sgM X} \vs$ for all vector fields $X$ and sections $\vs$ of $\vE$.
In other words $D^{(1)}$ obtained by restricting the invariant derivative to
tangent vectors using the isomorphism $\sgM\colon TM\dsum \pM\dsum T\dual M\to
\gM$ induced by $\gem$.
\item The \emphdef{Weyl connection} $D\colon \Cinf(M,\vE)\to \Cinf(M, T\dual
M\otimes \vE)$ is $D \vs = \gem^{-1}_\vE D^{(1)} (\gem_\vE \vs)$, i.e., the
connection on $\vE$ induced by $D^{(1)}$ on $\gr V$ via the isomorphism
$\gem_\vE\colon \vE\to \gr\vE$.
\end{numlist}
\end{defn}
By definition, $D^{(1)}$ and $D$ agree on bundles associated to semisimple
$P$-modules (when $\exp\p^\perp$ acts trivially and $\vE$ and $\gr\vE$ are
canonically isomorphic). In the notation we suppress their dependence on the
Weyl structure $\gem$.  A priori they also depend on the chosen $P$-module
$\aE$. This latter dependence is straightforward as they are associated to
principal $P$-connections. To see this, we use the isomorphism $\sgG\colon
\cG\cross(\T\dsum \p_0\dsum \T\dual)\to\cG\cross\g$ defined by the Weyl
structure to decompose the Cartan connection $\cc\colon T\cG\to\g$ as
\begin{equation}
\cc=\sgG\cc_\T + \cc_\p,\quad \cc_\p= \sgG\cc_{\p_0}+\cc_{\T\dual},
\end{equation}
where $\cc_\T:=(\cc \mod\p)\colon T\cG\to \T$ is the solder form, induced by
projecting $\cc$ onto $\T=\g/\p$ and similarly
$\cc_{\p_0}:=(\cc_{\p}\mod\p^\perp)\colon T\cG\to\p_0=\p/\p^\perp$. Thus
\begin{equation}
\sgG^{-1}\circ\cc=\cc_\T +\cc_{\p_0}+\cc_{\T\dual}
\end{equation}
and we can write the $\p$-part conceptually as
\begin{equation}
\cc_{\p}=(\sgG\cc_{\p_0}-\Ge^*\MCw)+(\cc_{\T\dual}+\Ge^*\MCw),
\end{equation}
where $\MCw$ is the Maurer-Cartan form on $\weyl$.  This leads to the
following proposition.

\begin{prop}\label{princ-pov}
Let $(\cG\to M,\cc)$ be a parabolic geometry with Weyl structure $\gem$.
Then\textup:
\begin{numlist}
\item $\cc_{\p}$ is a principal $P$-connection on $\cG$ inducing
$D^{(1)}$ on associated bundles\textup;
\item $\cc_{\Ge}=\sgG\cc_{\p_0}-\Ge^*\MCw$ is a principal $P$-connection on
$\cG$ inducing $D$ on associated bundles\textup;
\item $\rho=\cc_{\T\dual}+\Ge^*\MCw$ is a horizontal, $P$-invariant
$\p^\perp$-valued $1$-form on $\cG$\textup; if $r^D$ is the induced $T\dual
M$-valued $1$-form on $M$ and $\act$ is the natural action of $T\dual
M=\pM^\perp$ on $\vE$, then
\begin{equation}
D^{(1)}_X \vs = D\low_X \vs + r^D(X)\act \vs.
\end{equation}
\end{numlist}
\end{prop}
\begin{numlproof}
\item Clearly $\cc_{\p}$ is $P$-invariant and $\p$-valued, and so, since $\cc$
is a Cartan connection, $\cc_{\p}$ is a principal $P$-connection. Let $X$ be a
vector field and $\vs$ a section of $\vE$, and let $\chi\colon\cG\to \T$ and
$\vf\colon \cG\to \aE$ be the corresponding $P$-invariant functions. Since the
identification of $TM$ with $\cG\cross_P \T$ is via the solder form,
$\chi=\cc_\T(\hat X)$ for any $P$-invariant lift $\hat X$ of $X$ to $\cG$.

As a $P$-invariant function on $\cG$, $D^{(1)}_X\vs$ is then
\begin{equation*}
\InvDer_{\sgG \chi}\vf= \InvDer_{\sgG\cc_\T(\hat X)}\vf
= \InvDer_{\cc(\hat X)}\vf-\InvDer_{\cc_\p(\hat X)}\vf
= d\vf (\hat X)+ \cc_\p(\hat X)\act \vf,
\end{equation*}
which is precisely the $P$-invariant function on $\cG$ corresponding to the
covariant derivative of $\vs$ along $X$ induced by $\cc_\p$.

\item
We now mirror the construction of $D$ from $D^{(1)}$ on the principal bundle
level, using Remark~\ref{keyrem}\ref{kr2}: if $\aE$ is a filtered $P$-module
and $\vf\colon\cG\to \aE$ is $P$-invariant, corresponding to a section $\vs$
of $\vE=\cG\cross_P \aE$, then $D\vs$ corresponds to the $P$-invariant
horizontal $1$-form
\begin{equation*}
\Ge^{-1}_\aE (d+\cc_\p) (\Ge_\aE \vf)
=d\vf+\sgG(\cc_{\p_0})\act\vf+\Ge^{-1}_{\aE}(d\Ge_\aE)\vf
=d\vf+\cc_\Ge\act \vf
\end{equation*}
with $\cc_\Ge = \sgG(\cc_{\p_0})-\Ge^*\MCw$, by Lemma~\ref{MC}, as required.

\item This follows immediately because $\cc_\p=\cc_\Ge+\rho$. (One can also
easily see directly that $\cc_{\T\dual}+\Ge^*\MCw$ is a $P$-invariant
$\p^\perp$-valued horizontal $1$-form on $\cG$.)
\end{numlproof}

In conformal geometry $r^D$ is the \emph{normalized Ricci curvature} (aka.\
the \emph{Schouten} or \emph{Rho tensor}) of $D$. This is the origin of the
term Ricci-corrected Weyl connection.

We wish to see how the objects we have constructed depend on the choice of
Weyl structure. We can either do this on $M$, or for the corresponding
$P$-invariant objects on $\cG$.
\begin{prop}\label{change} Let $\gem$ and $\Tgem=(\Ad q^{-1})\gem$ be Weyl
structures, with $q\colon M \to \cG\cross_P \exp\p^\perp$ \textup(associated
to the adjoint action\textup), and let $\vE$ be a filtered
$P$-bundle. Then $\Tgem_\vE = \gem_\vE\circ (q\act)$.
\end{prop}
This is immediate from equation~\eqref{gr-change}. In practice it suffices to
understand infinitesimal variations.  Let $q_t$ be a curve of sections of
$\pM$ with $q_0=\iden$ and $\dot q_0=\ag$ for a $1$-form $\ag$ (equivalently a
$P$-invariant function $\cG\to\p^\perp$).  Then for any object $F(\gem)$
depending on $\gem$, define $(\pd_\ag F)(\gem)$ to be the $t$-derivative of
$F((\Ad q_t^{-1})\gem)$ at $t=0$ (so that $\pd_\ag\gem=-\ag$). By the
fundamental theorem of calculus, $F$ is independent of the Weyl structure if
and only if $(\pd_\ag F)(\gem)=0$ for all Weyl structures $\gem$ and all
$1$-forms $\ag$.  Proposition~\ref{change} implies that
$\pd_\ag\gem_\vE=\gem_\vE\circ(\ag\act)$. It is now straightforward to
differentiate the definition of $D^{(1)}$.
\begin{prop} For a Weyl structure $\gem$, a $1$-form $\ag$, and a
vector field $X$, let $[\ag,X]^{\gem}_{\pM}=[\ag,\sgM X]-\sgM(\ag\act X)$.
Then for any section $\vs$ of a filtered $P$-module $\vE$,
\begin{equation}\label{D1-change}
\pd_\ag D^{(1)}_X\vs = [\ag,X]^{\gem}_{\pM}\act \vs
\end{equation}
\end{prop}
\begin{proof} For $X\in TM$, $\pd_\ag\sgM X= \sgM(\ag\act X)-[\ag,\sgM X]=
-[\ag,X]^{\gem}_{\pM}$. This is $\pM$-valued (it is the $\pM$ component of the
Lie bracket $[\ag,X]$, where the lift of $X$ to $\gM$ and the projection to
$\pM$ are defined using $\gem$), and so $\pd_\ag \InvDer_{\sgM X}\vs
=\InvDer_{\pd_\ag\sgM X}\vs= [\ag,X]^{\gem}_{\pM}\act\vs$.
\end{proof}

The Ricci-corrected first derivative $D^{(1)}$ agrees with the Weyl connection
on semisimple modules. Ricci corrections start to play a more important role
when higher derivatives are considered, because jet modules are not
semisimple. The following deceptively simple definition clearly generalizes
the previous definition of $D^{(1)}$.

\begin{defn} Let $\gem$ be a Weyl structure. Then we define an isomorphism
from
\begin{equation*}
\hat J^k_\cc\vE\leq \Dsum_{j=0}^k(\tens^j\gM^*)\tens\vE
\qquad\text{to}\qquad
\hat J^k_{\gem}\vE:=\Dsum_{j=0}^k(\tens^jT\dual M)\tens\vE
\end{equation*}
by sending $\phi=(\phi_0,\phi_1,\ldots\phi_k)$ to
$\psi=(\psi_0,\psi_1,\ldots\psi_k)$ where
$\psi_j=(\phi_j\circ\sgM)\restr{\tens^j TM}$.

If $\vs\in\Cinf(M,\vE)$, then the section of $\hat J^k_{\gem}\vE$ corresponding
to $\hj^k_\cc\vs$ is denoted $\hj^k_D\vs=(\vs,D^{(1)}\vs,D^{(2)}\vs,\ldots
D^{(k)}\vs)$. We call $\hj^k_D$ the \emphdef{Weyl jet operator}, and its
components $D^{(j)}$ the \emphdef{Ricci-corrected powers of the Weyl
connection}.
\end{defn}
An alternative description of the Weyl jet operator is obtained from the
obvious natural isomorphism between $\gr\hat J^k_\cc\vE$ and $\gr\hat
J^k_{\gem}\vE$. Then $\hj^k_D=\gem^{-1}_{\smash[b]{\hat
J^k_{\gem}\vE}}\gem_{\smash[b]{\hat J^k_\cc\vE}} \hj^k_\cc$. It follows that
\begin{equation}\label{jet-change}
\pd_\ag\hj^k_D\vs=\gem^{-1}_{\hat J^k_{\gem}\vE}\gem_{\hat J^k_\cc\vE} \ag\cdot
\hj^k_\cc\vs - \ag\cdot \gem^{-1}_{\hat J^k_{\gem}\vE}\gem_{\hat
J^k_\cc\vE}\hj^k_\cc\vs= \ag\ast\hj^k_D \vs - \ag\act\hj^k_D\vs,
\end{equation}
where $\ast$ is the action of $T\dual M$ on $\hat J^k_{\gem}\vE$ induced by
the isomorphism $\gem^{-1}_{\hat J^k_{\gem}\vE}\gem_{\hat J^k_\cc\vE}$ with
$\hat J^k_\cc\vE$.

The computation of $\ast$ is a straightforward exercise in algebra: it
suffices to describe the jet $\tilde\phi=\ag\act\phi$ for $\phi\in\hat
J^k_0\aE$, $\ag\in\p^\perp$, in terms of the elements $\tilde\psi$ and $\psi$
of $\hat J^k_{\age}\aE:=\Dsum_{j=0}^k(\tens^j\T\dual)\tens\aE$ corresponding
to $\tilde\phi$ and $\phi$ using an algebraic Weyl structure $\age$. This
$\p$-module structure on $\hat J^k_{\age}\aE$ is computed
in~\cite{CSl,Slovak}: let us sketch briefly the computation. We write
$\psi=(\psi_0,\psi_1,\ldots\psi_k)$ with $\psi_j\in (\tens^j\T\dual)\tens\aE$
and similarly for $\tilde\psi$.

First, note that the existence of natural projections $\hat J^k_0\aE\to \hat
J^\ell_0\aE$ (for $k\geq\ell$) implies that the $\psi_j$ component contributes
only to $\tilde\psi_{j+s}$ for $s\geq0$, independently of $k\geq j+s$. We may
therefore write $\tilde\psi_\ell=(\ag\ast\psi)_\ell
=\sum_{j+s=\ell}\ag\ast_s\psi_j$ for any $\ell\leq k$, where $\ag\ast_s\psi_j$
denotes the contribution from $\psi_j$. Clearly
$\ag\ast_0\psi_j=\ag\act\psi_j$ is the ordinary action of $\p^\perp\leq\p$.

Secondly, observe that the natural inclusions $\hat J^{k+1}_0\aE\to J^1_0\hat
J^k_0\aE$ mean that we can compute $\ag\ast_s\psi_j$ inductively using the
definition of $J^1_0\aE$. Recall that this consists of the pairs
$(\phi_0,\phi_1)$ in $\aE\dsum(\g\dual\tens\aE)$ with
$\phi_1(\xi)+\xi\act\phi_0=0$ for $\xi\in\p$.  The identification with
$\aE\dsum(\T\dual\tens\aE)$ is obtained by restricting $\phi_1$ to $\T$, using
$\age$. The induced $\p$-module structure on $\aE\dsum(\T\dual\tens\aE)$ is
given by $\xi\act(\phi_0,\phi_1\restr{\T})=\bigl(\xi\act\phi_0,
(\xi\act\phi_1)\restr{\T}\bigr)$, and one easily computes that
$(\xi\act\phi_1)\restr{\T}=\xi\act(\phi_1\restr{\T})
+[\xi,\cdot]_\p^{\age}\act\phi_0$, where $\age$ is used to split the natural
maps $\p\to\g\to\T$. Therefore
\begin{equation*}
\ag\ast(\psi_0,\psi_1)=
(\ag\act\psi_0,\ag\act\psi_1+\ag\ast_1\psi_0)\qquad\text{with}\qquad
\ag\ast_1\psi_0= [\ag,\cdot]_{\p}^{\age}\act\psi_0.
\end{equation*}
(Applying this pointwise on $M$, we rederive equation~\eqref{D1-change}.)
Iterating this action we have
$\ag\ast_s\psi_0=\bigl[[\ldots[\ag,\cdot]_{\p}^{\age},
\ldots]_{\p}^{\age},\cdot\bigl]{}_{\p}^{\age}\act\psi_0$.  The formula for
$\ag\ast_s\psi_j$ is also computed inductively, by considering the action of
$\ag$ on $J^1_0\hat J^{j+s-1}_0\aE$: apply this action to
$(0,\ldots0,\phi_j,0,\ldots 0)\in\hat J^{j+s-1}_0\aE$ to obtain
$\ag\ast_s\psi_j$ as the $\T\dual\tens(\tens^{j+s-1}\T\dual)\tens\aE$
component with respect to $\age$.

Passing from the algebra to associated bundles, we obtain the following
result.
\begin{prop}\label{jetact}
The action of $\ag\in T\dual M$ on $\psi\in\hat J^k_{\gem} \vE$
induced by the identification with $\hat J^k_\cc \vE$ is given by
\begin{equation*}
(\ag\ast\psi)_\ell=\sum_{j+s=\ell}\ag\ast_s\psi_j
\end{equation*}
where, if we suppose that $\psi_j=\td_1\tens\cdots\tens\td_j\tens\ve$ for
$\td_i\in T\dual M$ and $\ve\in \vE$, we have
\begin{align*}
\ag\ast_0\psi_j&=\ag\act\psi_j\\
\ag\ast_1\psi_j&=\sum_{0\leq i\leq j}\td_1\tens\td_2\tens
\cdots\tens\td_{i}\tens[\ag,\cdot]_{\pM}^{\gem}\act(\td_{i+1}\tens\cdots
\tens\td_j\tens\ve)\\
\ag\ast_s\psi_j&=\smash{\sum_{\substack
{0\leq i_1\leq i_2\leq\cdots\leq i_s\leq j\\ a_1,a_2,\ldots a_s}}}\;\,
\td_1\tens\cdots\tens\td_{i_1}\tens\dfr^{a_1}\tens\td_{i_1+1}\tens\cdots
\tens\td_{i_2}\tens\dfr^{a_2}\tens\td_{i_2+1}\tens\cdots\\
&\qquad\qquad\qquad\qquad
\cdots\tens\cdots\tens\td_{i_{s-1}}\tens\dfr^{a_{s-1}}
\tens\td_{i_{s-1}+1}\tens\cdots\tens\td_{i_s}\tens\dfr^{a_s}\tens\\
&\qquad\qquad\qquad
\bigl[[\ldots[[\ag,\frm_{a_1}]_{\pM}^{\gem},\frm_{a_2}]_{\pM}^{\gem},\ldots
\frm_{a_{s-1}}]_{\pM}^{\gem},\frm_{a_s}\bigr]{}_{\pM}^{\gem}
\act(\td_{i_s+1}\tens\cdots\tens\td_j\tens\ve).
\end{align*}
In these formulae $\frm_a$ is a frame of $TM$ with dual frame $\dfr^a$.
\end{prop}
This action not only gives an explicit formula for $\pd_\ag \hj_D$: it also
provides an explicit inductive formula for the Ricci-corrected powers
$D^{(k)}$ of the Weyl connection $D$. Indeed, since the order $k+1$ part of
$j^1_\cc\hj^k_\cc\vs$ is the same as that of $\hj^{k+1}_\cc \vs$, we obtain:
\begin{equation}\label{rcdexp}\begin{split}
\iota\low_X D^{(k+1)}\vs&=D\low_X D^{(k)}\vs
+\project_{(\tens^k T\dual M)\tens\vE}\bigl(r^D(X)\ast\hj^k_D\vs\bigr)\\
&=D\low_X D^{(k)}\vs+\sum_{j+s=k} r^D(X)\ast_s D^{(j)}\vs.
\end{split}\end{equation}
Explicit formulae for invariant differential operators will follow by
computing projections of $D^{(k)}$ using this inductive expression and some
representation theory.

\section{Strongly invariant operators}\label{sio}

\subsection{Jet module homomorphisms}

Our goal is to explain how Ricci-corrected Weyl differentiation leads to
explicit formulae for a class of invariant differential operators.  These are
the \emphdef{strongly invariant} operators of~\cite{CSS3,CSS4,ES}, defined as
follows.  Let $\aE$ and $\aF$ be $P$-modules and let $\Phi\colon\hat
J^k_0\aE\to\aF$ be a $P$-homomorphism.  Then $\Phi$ induces a bundle map
$F\colon\cG\cross_P\hat J^k_0\aE\to\cG\cross_P\aF$, and hence an invariant
differential operator $F\circ\hj^k_\cc$ from $\vE$ to $\vF$.

In practice, such $P$-homomorphisms are constructed by lifting a
$P_0$-homomorphism $\gr\Phi\colon\gr\hat J^k_0\aE\to \gr \aF$ using an
algebraic Weyl structure $\age$ to give
$\Phi_{\age}=\age^{-1}_\aF\circ\gr\Phi\circ\age_{\aE}$.  Since $\Phi_{(\Ad
q)\age}\ve=q\act \Phi_{\age}(q^{-1}\act\ve)$, it follows that
$\Phi=\Phi_{\age}$ is a $P$-homomorphism if and only if it is independent of
the algebraic Weyl structure $\age$.

If we mirror this construction on associated bundles, for any Weyl structure
$\gem$, a bundle map $\gr F\colon \gr\hat J^k\vE\to \gr\vF$ (associated to
$\gr\Phi\colon \gr \hat J^k_0\aE\to \gr \aF$) induces a differential
operator $\gem^{-1}_\vF\circ\gr F\circ\gem_{\hat J^k\vE}\circ\hj^k_\cc$
from $\vE$ to $\vF$, which will be invariant if it is independent of the
choice of Weyl structure $\gem$. An obvious sufficient condition is that
$\gem^{-1}_\vF\circ\gr F\circ\gem_{\hat J^k\vE}$ is independent of the
choice of Weyl structure and these are the strongly invariant operators. (The
condition is not necessary because $\hj^k_\cc\vs$ will satisfy some Bianchi
identities not satisfied by general sections of $\hat J^k\vE$.)

The method we shall adopt for constructing strongly invariant operators is to
construct a $P$-homomorphism $\Phi\colon \hat J^k_{\age}\aE\to \aF$, inducing
a bundle map $F\colon \hat J^k_{\gem}\vE\to \vF$ and hence, for any Weyl
structure $\gem$, a differential operator $F\circ\hj^k_D$ from $\vE$ to
$\vF$. Since $\hat J^k_{\age}\aE$ has the same associated graded module as
$\hat J^k_0\aE$, it is straightforward to obtain conditions that $\Phi$
induces a $P$-homomorphism $\hat J^k_0\aE\to \aF$, and hence a strongly
invariant operator.  The expression $F\circ\hj^k_D$ then gives an explicit
formula for this operator in terms of a Weyl structure.

The $P$-homomorphisms we construct here all factor through the projections
$\hat J^k_{\age}\aE\to \tens^k \T\dual\tens \aE\to \tens^k\f\dual\tens\aE$,
where the second projection is induced by the restriction map
$\T\dual\to\f\dual$. Our first task is to apply these projections to the jet
module action $\ast$. Then we apply further projections
$(\tens^k\f\dual)\tens\aE\to\aF$ to obtain the $P$-homomorphisms we seek. In
geometric terms, operators obtained from homomorphisms factoring through these
projections are given by applying a bundle map $(\tens^k\fM^*)\tens\vE\to\vF$
to the restriction of $D^{(k)}$ to $\tens^k\fM$.

\subsection{Restricting the jet module action}\label{pjm}

We first restrict $\ag\ast\psi$ to $\f$. For notational simplicity we give the
result algebraically: the formulae on associated bundles easily follow.

\begin{prop} Suppose $\aE$ is a semisimple $P$-module. Let $\rtf$
denote the restriction maps $\tens^j\T\dual\to\tens^j\f\dual$ and let
$\frm_a$, $\dfr^a$ be dual bases for $\f$, $\f\dual$.  Then
$(\ag\ast\psi)_\ell=\sum_{j+s=\ell}\ag\ast_s\psi_j$, where for each $j$ we
have, supposing $\psi_j=\td_1\tens\cdots\tens\td_j\tens\ve$ as before,
\begin{align*}
\rtf(\ag\ast_0\psi_j)&=0\\
\rtf(\ag\ast_1\psi_j)&=\sum_{0\leq i\leq j}\rtf(\td_1\tens\td_2\tens
\cdots\tens\td_{i})\tens\Psi\bigl(\ag\tens\rtf(\td_{i+1}\tens\cdots
\tens\td_j\tens\ve)\bigr)\\
\rtf(\ag\ast_s\psi_j)&=\smash{\sum_{\substack
{0\leq i_1\leq i_2\leq\cdots\leq i_s\leq j\\ a_1,a_2,\ldots a_{s-1}}}}\;\,
\rtf(\td_1\tens\cdots\tens\td_{i_1})\tens\dfr^{a_1}\tens
\rtf(\td_{i_1+1}\tens\cdots\tens\td_{i_2})\tens\dfr^{a_2}\tens
\cdots\\
&\qquad\qquad\qquad\qquad\qquad
\cdots\tens\cdots\tens\dfr^{a_{s-1}}\tens
\rtf(\td_{i_{s-1}+1}\tens\cdots\tens\td_{i_s})\tens\\
&\qquad\qquad\qquad
\Psi\bigl(\bigl[\ldots[[\ag,\frm_{a_1}]_\p^{\age},\frm_{a_2}]_\p^{\age},\ldots
\frm_{a_{s-1}}\bigr]{}_{\f\dual}^{\age}\tens\rtf(\td_{i_s+1}\tens\cdots
\tens\td_j
\tens\ve)\bigr).
\end{align*}
Here, for any semisimple $P$-module $\tilde\aE$, we define
$\Psi\colon\f\dual\tens\tilde\aE\to\f\dual\tens\tilde\aE$ by
\begin{equation}
\Psi(\td\tens\tilde\ve)(\chi)=[\td,\chi]\act\tilde\ve
\end{equation}
for $\chi\in\f$, which is well defined since $\tilde\aE$ is semisimple. In the
above formulae for $\ast_1$ and $\ast_s$, we have
$\tilde\aE=(\tens^{j-i}\f\dual)\tens\aE$ and $(\tens^{j-i_s}\f\dual)\tens\aE$
respectively.
\end{prop}
\begin{proof} This is immediate from (the algebraic version of)
Proposition~\ref{jetact}, and the equality
\begin{multline*}
\rtf\bigl(\bigl[[\ldots[[\ag,\frm_{a_1}]_\p^{\age},\frm_{a_2}]_\p^{\age},\ldots
\frm_{a_{s-1}}]_\p^{\age},\frm_{a_s}\bigr]{}_\p^{\age}
\act(\td_{i_s+1}\tens\cdots\tens\td_j\tens\ve)\bigr)\\
=\bigl[[\ldots[[\ag,\frm_{a_1}]_\p^{\age},\frm_{a_2}]_\p^{\age},\ldots
\frm_{a_{s-1}}]_{\f\dual}^{\age},\frm_{a_s}\bigr]{}_{\p_0}^{\age}\act
\bigl(\rtf(\td_{i_s+1}\tens\cdots\tens\td_j)\tens\ve\bigr),
\end{multline*}
which holds because the projection of the Lie bracket
$\p^\perp\tens\p^\perp\to\p^\perp\to\f\dual$ vanishes and the module
$(\tens^{j-i_s}\f\dual)\tens\aE$ is semisimple.
\end{proof}
The formulae of this proposition may not seem simpler than the full formulae
of Proposition~\ref{jetact}, but they are easier to handle, since $\Psi$ is an
operator on semisimple modules.

\subsection{Special types of projections}

To progress further, we need to understand the map
\begin{equation*}
\Psi\colon\f\dual\tens\aE\to\f\dual\tens\aE,\quad
\Psi(\td\tens\ve)=\sum_{a}\dfr^a\tens([\td,\frm_a]\act\ve),
\end{equation*}
where $\aE$ is a semisimple $P$-module and $\frm_a$ and $\dfr^a$ are dual
bases of $\f$ and $\f\dual$ respectively.  Since $\Psi$ is a $P$-homomorphism,
it acts by a scalar on every irreducible component of $\f\dual\tens\aE$ and
these scalars can be computed explicitly using Casimirs.  Note now that all
weights of the $\p_0$-module $\T\dual$ have multiplicity one (they are just
positive roots of $\g$). Hence results from~\cite{Bourbaki,PRV,SlSo} show that
all irreducible components of $\T\dual\tens\aE$ are multiplicity free.

We can write $\f\dual\tens\aE$ as a sum of well-defined irreducible
$P_0$-components $\aE_{\ind}$.  Consequently,
$(\tens^2\f\dual)\tens\aE=\f\dual\tens(\Dsum_{\ind_1}\aE_{\ind_1})$ can be
again written as a sum of invariant subspaces labelled by a couple
$(\ind_1,\ind_2)$ of indices indicating that it is the isotypic component with
label $\ind_2$ inside $\f\dual\tens\aE_{\ind_1}$. Inductively, we get a well
defined decomposition of $(\tens^k\f\dual)\tens\aE$ into $P_0$-invariant
subspaces labelled by paths $\ind=(\ind_1,\ldots\ind_k)$ of indices showing
their positions in consecutive decompositions.  Note that the full isotypic
component of $(\tens^k\f\dual)\tens\aE$ with label $\ind_k$ is the direct sum
of all subspaces $\aE_{\ind}$ labelled by paths $\ind$ ending with $\ind_k$.

We now suppose that our component $\aE_{\ind}$ is in the \emph{symmetric}
tensor product $S^k\f\dual\tens\aE$.

\begin{prop}\label{special}
Let $\pi$ be a projection of $\hat J^k_{\age}\aE$ to an invariant subspace in
$(S^k\f\dual\tens\aE)\intersect\aE_{\ind}$ for some $\aE_{\ind}$ in the
decomposition described above. Then for any element $\psi$ of $\hat
J^k_\eps\aE$, the only contribution to $\pi(\ag\ast\psi)$ is from $\psi_{k-1}$
and if $\psi=\td_1\tens\cdots\tens\td_{k-1}\tens\ve$ we have
\begin{equation*}
\pi(\ag\ast\psi)=\sum_{0\leq i\leq k-1}\pi\bigl(\td_1\tens\td_2\tens
\cdots\tens\td_{i}\tens\Psi\bigl(\rtf(\ag)\tens\rtf(\td_{i+1}\tens\cdots
\tens\td_{k-1}\tens\ve)\bigr)\bigr).
\end{equation*}
\proofof{prop} Note first that $\pi(\ag\ast\psi)=\sum_{j+s=k}
\pi(\ag\ast_s\psi_j)$ and that the $s=0$ term vanishes (since $\aE$ is
semisimple). Hence it remains to show that the terms with $s\geq2$ are zero.
To do this is suffices to show that terms in which $\Psi$ is applied to a Lie
bracket are killed by the projection. Consider therefore an expression of the
form
\begin{equation*}
\textstyle\sum_a\dfr^a\tens\Psi([\ag,\frm_a]_{\f\dual}^{\age}\tens\ve)
\end{equation*}
and suppose we apply a projection $\pi$ to $(\tens^2\f\dual)\tens\aE$ which
factors through $S^2\f\dual\tens\aE$ and is of the form
$\pi_2\circ\iden\tens\pi_1$ where $\pi_1\colon\f\dual\tens\aE\to\aE_1$ is a
projection onto an isotypic component. Since $\Psi$ acts by a scalar on such
a component, the result is a multiple of $\pi$ applied to
\begin{equation*}
\textstyle\sum_a\dfr^a\tens [\ag,\frm_a]_{\f\dual}^{\age}\tens\ve.
\end{equation*}

The projection $\pi$ factorizes through the symmetric product, so it
suffices to note that
\begin{equation*}
\Bigl({\textstyle\sum_{a}\dfr^a\tens[\ag,\frm_a]_{\f\dual}^{\age}}\Bigr)
(\chi_1,\chi_2)=
[\ag,\chi_1]_{\f\dual}^{\age}(\chi_2)=\ag([\sg\chi_1,\sg\chi_2]\mod\p),
\end{equation*}
which is clearly antisymmetric in $\chi_1,\chi_2\in\f$. The terms appearing in
the action for $s\geq2$ are all of this form, with $\ag$ and $\ve$ replaced by
iterated brackets and suitable tensor products.
\end{proof}

\subsection{Casimir computations}\label{cc} In this section we compute the
eigenvalues of $\Psi$ and hence prove that $\Psi$ acts by a scalar on each
isotypic component. This is the first point at which we need detailed
information from representation theory, so we set up the necessary notation.

Let $\g$ be a semisimple Lie algebra with complexification $\g^{\C}$.  We
choose a Cartan subalgebra $\h^\C$ in $\g^{\C}$, a set $\Delta^+$ of positive
roots, and its subset $S=\{\alpha_1,\ldots\alpha_r\}$ of simple roots.  Using
the Killing form $(\cdot\,,\cdot)$, with any normalization, the fundamental
weights $\omega_1,\ldots\omega_r$ are defined by
$(\alpha_i^\vee,\omega_j)=\de_{ij}$, where $\alpha_i^\vee
=2\alpha_i/(\alpha_i,\alpha_i)$.

The dominant Weyl chamber $\cC$ is given by linear combinations of fundamental
weights with nonnegative coefficients.  Finite dimensional complex irreducible
representations of $\g^{\C}$ (as well as of $\g$) are characterized by their
highest weights $\lam\in\cC$, which lie in the weight lattice
$\{\sum\lam_i\omega_i: \lam_i\in\Z\}$.  The corresponding representation will
be denoted by $\aE_{\lam}$.

A reductive algebra is a direct sum of a commutative and a semisimple algebra
(either of which can be trivial).  Its irreducible (complex) representations
are tensor products of irreducible representations of the summands, where
irreducible representations of a commutative Lie algebra $\lie{a}$ are one
dimensional, characterized by an element of Hom$(\lie{a},\C)$.

\begin{rem} For simplicity, we focus on complex representations of the
Lie algebras in question.  In practice, we may well be more interested in real
representations. For this, it is sufficient to use the following description:
a real or quaternionic structure on a complex $\g$-module $\aE$ is a
conjugate-linear $\g$-map $J\colon\aE\to \aE$ with $J^2=1$ or $J^2=-1$
respectively. An irreducible real representation of a Lie algebra can be
identified either with an irreducible complex representations, or with such a
representation endowed with a real or quaternionic structure. We note also
that if $\V$ is a complex $\g$-module and $\U$ is a real $\g$-module with
complexification $\U^c$, then $\U\otimes\V$ and $\U^c\otimes_\C \V$ are
equivalent as complex modules.
\end{rem}

We now specialize to the situation where $\g$ is a semisimple Lie algebra with
parabolic subalgebra $\p$ and Levi factor $\p_0$.  The set $S$ of simple roots
for $\g^{\C}$ can be chosen in such a way that all positive root spaces are
contained in $\p^{\C}$. This fixes an algebraic Weyl structure, and the
positive root spaces lying in $\p^{\C}_0$ correspond to roots in the span of
the subset $S_0$ of `uncrossed' simple roots---we write $S=S_\times\cup S_0$
for the decomposition into crossed and uncrossed simple roots.  In this
situation, we shall say that a weight $\lam$ (integral for $\g$) is dominant
for $\p$, if its restriction to $\h^\C_{ss}=\h^\C\cap\p^{\C}_{0,ss}$ is
dominant for $\p^{\C}_{0,ss}$. Such a weight specifies uniquely an irreducible
$P$-module.

Let us denote by $\de$ the half sum of positive roots for $\g^{\C}$ and by
$\de_0$ the half sum of those positive roots for $\g^{\C}$ for which the
corresponding root space belongs to $\g^{\C}_0$.

\begin{prop}\label{casimir}
Let $\aE_{\lam}$ be an irreducible representation of $\p_0$ with highest
weight $\lam\in(\h^{\C})^*$. Let $\f\dual\tens\aE_{\lam}=\Dsum_{\mu\in
A}\aE_{\mu}$ be the decomposition of the tensor product into the sum of
isotypic components with highest weight $\mu$ and let $\pi_{\mu}$ be the
projection to $\aE_{\mu}$.  Let $\de$ denote the half sum of positive roots
for the Lie algebra $\g$. Then
\begin{align*}
\Psi&=\sum c_{\mu}\pi_{\mu},\\
\tag*{\textit{with}}
c_{\mu}&=\tfrac{1}{2}\left(|\mu+\de|^2-|\lam+\de|^2\right)
\end{align*}
and $|\alpha|^2=(\alpha,\alpha)$. \textup[Note that $\Psi$ depends upon the
normalization of the Killing form, since we used it to identify $\p^\perp$
with $(\g/\p)\dual$, and hence the bracket of $\f\dual$ with $\f$ depends on
$(\cdot\,,\cdot)$.\textup]

\proofof{prop} We first give a formula for $\Psi$ in terms of the Casimir
operator $C$ of $\p_0$, as in~\cite{Fegan,SlSo}. Let $E^i$ and $E_i$ be
bases for $\p_0$ which are dual with respect to $(\cdot\,,\cdot)$. Then
\begin{align*}
\Psi(\ag\tens\ve)&={\textstyle\sum_a \dfr^a\tens[\ag,\frm_a]\act\ve}
={\textstyle\sum_i [E^i,\ag]\tens E_i\act\ve}
={\textstyle\sum_i [E_i,\ag]\tens E^i\act\ve}\\
\tag*{and so}
2\Psi(\ag\tens\ve)&={\textstyle\sum_i E^i\act E_i\act(\ag\tens\ve)}
-\ag\tens\Bigl({\textstyle\sum_iE^i\act E_i\act\ve}\Bigr)-
\Bigl({\textstyle\sum_iE^i\act E_i\act\ag}\Bigr)\tens\ve\\
&=C(\ag\tens\ve)-\ag\tens C(\ve)-C(\ag)\tens\ve.
\end{align*}
It is well known that on an irreducible representation with highest weight
$\lam$, $C$ acts by the scalar $(\lam,\lam+2\de_0)$: with the definition of
$\de_0$ above, this holds even though $\p_0$ is reductive rather than
semisimple---see~\cite{SlSo}.

Let us write $S_\times=\{\be_i:i=1,\ldots r_\times\leq r\}$ for the crossed
simple roots of $\g$. We know that $-\be_i$ are precisely the highest weights
of the irreducible components $\f\dual_i$ of the $P_0$-module $\f\dual$, so
that the irreducible components of $\f\dual_i\tens\aE_{\lam}$ have highest
weights of the form $\lam-\be_i+\gam$, where $\gam$ is an integral linear
combination of simple roots for $\p^{\C}_0$. Hence for any isotypic component
$\aE_{\mu}$ there is an $i$ so that $\aE_\mu$ is an invariant subspace of
$\f\dual_i\tens\aE_{\lam}$.

Since the action of the Casimir depends only on the highest weight, we deduce,
following~\cite{SlSo}, that $\Psi$ acts on the entire isotypic component
$\aE_{\mu}$ by the scalar
\begin{equation*}
c=\tfrac12\bigl[(\mu,\mu+2\de_0)-(\lam,\lam+2\de_0)
-(-\be_i,-\be_i+2\de_0)\bigr].
\end{equation*}
It remains to identify $c$ with the constant $c_\mu$ above. For any
simple root $\be$, we have
\begin{equation*}
2(\de,\be)=(\de,2\be/|\be|^2)\,|\be|^2=
\sum_{j=1}^n(\omega_j,\be^\vee)|\be|^2=|\be|^2.
\end{equation*}
Hence $(\mu,\mu+2\de_0)-(\lam,\lam+2\de_0) -(-\be_i,-\be_i+2\de_0)$ is
given by
\begin{multline*}
(\mu,\mu+2\de)-(\lam,\lam+2\de)-2(\de-\de_0,\mu-\lam) -2(\be_i,\de-\de_0)\\
=|\mu+\de|^2-|\lam+\de|^2+2(\de-\de_0,-\be_i-\mu+\lam).
\end{multline*}
We know that $\mu-\lam$ is a weight of $\f\dual_i,$ hence
$-\be_i-(\mu-\lam)=\sum_{\alpha\in S_0}n_{\alpha}\alpha$. But for all
$\alpha\in S_0,$ we have $(\alpha,\de-\de_0)=0$, so the last term vanishes.
\end{proof}

\section{Explicit constructions of invariant operators}\label{ecio}

In this section, we are going to construct a large class of invariant
differential operators.  Most of them belong to the class of standard regular
operators, but a certain subclass are standard singular operators.  The class
of operators constructed here does not cover all standard regular operators,
but we shall see in examples that it covers many of them.

Various special cases of the main results of this section can be found
in~\cite{Branson,CSS3,CSS4,Fegan,SlSo}.
\begin{thm}\label{th1}
Let $\lroot\in(\h^\C)^*$ be a positive root with $\g_{\lroot}\subset \f\dual$.
In the case that $\g$ has roots of different lengths, we shall suppose that
$\lroot$ is a long root. Let $\lam,\mu$ be two integral dominant weights of
$\p_0$ with the property
\begin{equation*}
\mu+\de=\si_{\lroot}(\lam+\de)=\lam+\de-(\lam+\de,\lroot^\vee)\lroot.
\end{equation*}
Interchanging $\lam$ and $\mu$ if necessary, suppose that
$k:=-(\lam+\de,\lroot^\vee)$ is positive.
\begin{numlist}
\item There is a unique irreducible component $\aE_{\mu}$ with highest weight
$\mu$ in $(\tens^k\T\dual)\tens\aE_\lam$. Furthermore, $\aE_\mu$ belongs to
$S^k\f\dual \tens\aE_{\lam}$ and is of the form $\aE_{\ind}$, where
$\aE_{\ind_j}=\aE_{\lam+j\lroot}$ for $\ind=(\ind_1,\ldots\ind_k)$.

\item If $\pi\colon\hat J^k_{\age}\aE_{\lam}\to\aE_{\mu}$ is the corresponding
projection, then $\pi$ induces a $P$-homomorphism
$J^k_0\aE_{\lam}\to\aE_{\mu}$ and hence a strongly invariant differential
operator $\pi\circ D^{(k)}$ of order $k$ from sections of $\vE_{\lam}$ to
sections of $\vE_{\mu}$.
\end{numlist}
\end{thm}
\begin{numlproof}
\item To prove uniqueness, let us note first that all weights of $\T\dual$ are
positive roots.  Weights of $\tens^k\T\dual$ are hence sums of $k$ positive
roots.  The highest weight of any irreducible components of $(\tens^k\T\dual)
\tens\V_{\lambda}$ is of a form $\lambda +\beta,$ where $\beta$ is a weight of
$\tens^k\T\dual.$ The unicity claim is therefore true by the triangle
inequality, $\lroot$ being a long root.
                                                          
We now show that there is such a component $\aE_\mu$. By assumption, both
$\lam$ and $\mu=\lam+k\lroot$ are $P$-dominant, and $j\lroot$ is an extremal
weight of $\x^j \f\dual_i$ for all $j=1,\ldots k$.  The so-called
Parthasarathy--Ranga-Rao--Varadarajan conjecture (proved in \cite{Kumar})
states that if $\lam,\nu$ are highest weights of two irreducible (complex)
$\p_0$-modules $\aE_{\lam},\aE_{\nu}$, and if $\lroot$ is an extremal weight
of $\aE_{\nu}$, then an irreducible component $\aE_{\lam+\lroot}$ with
extremal weight $\lam+\lroot$ will appear with multiplicity at least one in
$\aE_{\lam}\tens_\C\aE_{\nu}$.  (In concrete cases, more elementary arguments
are available.) It follows that there is an irreducible component of
$(\x^j\f\dual_i)\tens\aE_{\lam}$ having $\lam+j\lroot$ as its highest
weight. Hence necessarily $\aE_{\mu}\subset (\x^k
\f\dual_i)\tens\aE_{\lam}\subset S^k\f\dual_i\tens\aE_{\lam}$.  Furthermore,
the same holds with $\lam$ replaced by $\lam+j'\lroot$ for all $j'=1,\ldots
k-1$. Hence, by uniqueness, the projection factors through
$(\x^j\f\dual_i)\tens\aE_{\lam}$.

\item Let us prove now that $\pi$ induces a $P$-homomorphism.  The action of
$\T\dual$ on $\aE_{\mu}$ is trivial, hence we must show that
$\pi(\ag\ast\psi)$ vanishes for any $\psi\in \hat J^k_{\age}\aE_\lam$ and any
$\ag\in\T\dual$.

Since $\pi$ is the projection to an irreducible piece of
$S^k\f\dual\tens\aE_{\lam}$ lying in a component of the form $\aE_\ind$, the
action is given by Proposition~\ref{special}.  Using the Casimir computation
of Proposition~\ref{casimir}, the projection of the action by $\ag$ on an
element $\psi\in\hat J^k\aE$ is given by
\begin{align*}
\pi(\ag\ast\psi)&=c\,\pi(\ag\tens\psi_{k-1})\\
\tag*{with}
2c&=\sum_{j=1}^{k}
\bigl(\,|\lam+j\lroot+\de|^2-|\lam+(j-1)\lroot+\de|^2\,\bigr)=
|\mu+\de|^2-|\lam+\de|^2,
\end{align*}
which is zero because $\mu+\de=\si_{\lroot}(\lam+\de)$ and $\si_\lroot$ is an
isometry.
\end{numlproof}

We turn now the formulae for these operators in terms of a Weyl
structure. Explicit formulae for the coefficients of various curvature terms
for standard operators were first found in the conformal
case~\cite{Baston3,Gover}, and later extended to the $|1|$-graded
case~\cite{Baston12,CSS3}. A very surprising fact was that the formulae were
quite universal and did not depend on the specific parabolic structure or on
the highest weights of the representations involved. The general structure of
coefficients described in~\cite{CSS3} was quite complicated. Here we notice
that the organization of the terms produced by the Ricci-corrected derivative
leads to much simpler coefficients. At the same time, the formulae are
extended from the $|1|$-graded case to the broad class of standard operators
in all parabolic geometries with no extra complications: the form of the
operator depends only on its order.

\begin{thm}\label{th2} Let a positive integer $k$ and a long root $\lroot$
with $\g_\lroot\subset \f\dual$ be given.  Define differential operators
$\cD_{k,j}$ \textup(of order $j=0,\ldots k$\textup), acting on sections of any
associated bundle, by the recurrence relation
\begin{equation}\label{recur}
\iota\low_X \,\cD_{k,j+1}=D\low_X\circ \cD_{k,j}+j(k-j)\Gam(X)\tens
\cD_{k,j-1}
\end{equation}
with $\cD_{k,0}=\iden,\,\cD_{k,1}=D$.  Here $D$ be the covariant derivative
given by the choice of the Weyl structure and $\Gam=-\frac12|\lroot|^2r^D$.
Let $\cD_k=\cD_{k,k}$. Then any invariant operator of order $k$ constructed in
Theorem~\textup{\ref{th1}}, mapping sections of $\vE_{\lam}$ to sections of
$\vE_{\mu}$, $\mu=\lam +k\lroot$ is given by $\pi\circ \cD_k$ where $\pi$ is
the projection onto $\vE_{\mu}$.

\proofof{thm} We know that the invariant operator is given by $\pi\circ
D^{(k)}$ and we have given a recurrence formula for $D^{(k)}$ in
section~\ref{rcd}.  Hence we only have to compute the projection of the action
of $r^D$ on $S^j\f\dual\tens \aE_{\lam}$, which is straightforward using the
results of \S\ref{pjm}-\ref{cc}: we find that the action is the tensor product
with $c\,r^D,$ where $2c=|\lam+j\lroot+\de|^2-|\lam+\de|^2$.

Now, since $k=-(\lam+\de,\lroot^\vee)$, we have
\begin{equation*}
2c=|\lam+j\lroot+\de|^2-|\lam+\de|^2=(2\lam+j\lroot+2\de,j\lroot)
=|\lroot|^2j(-k+j).
\end{equation*}
Substituting this into the projection of the recurrence formula for
$D^{(j+1)}$ gives~\eqref{recur}.
\end{proof}
Hence the universal nature of the explicit formulae for invariant operators
arises from the fact that $\cD_k$ only depends upon $k$.  It is
straightforward to compute $\cD_k$ for small $k$. Omitting the tensor product
sign when tensoring with $\Gam^j=\Gam\tens\cdots\tens\Gam$, we have
\begin{align*}
\cD_1s=Ds{}^{\hphantom{1}}&\\
\cD_2s=D^2s&+\Gam s\\
\cD_3s=D^3s&+2D(\Gam s)+2\Gam Ds\\
\cD_4s=D^4s&+3D^2(\Gam s)+4D(\Gam Ds)+3\Gam D^2s\\
&+9\Gam^2 s
\displaybreak[0]\\
\cD_5s=D^5s&+4D^3(\Gam s)+6D^2(\Gam Ds)+6D(\Gam D^2s)+4\Gam D^3s\\
&+24D(\Gam^2 s)+16\Gam D(\Gam s)+24\Gam^2 Ds\\
\cD_6s=D^6s&+5D^4(\Gam s)+8D^3(\Gam Ds)+9D^2(\Gam D^2s)+8D(\Gam D^3s)
+5\Gam D^4s\\
&+45D^2(\Gam^2 s)+40D(\Gam D(\Gam s))+25\Gam D^2(\Gam s)\\
&\qquad\qquad+64D(\Gam^2 Ds)+40\Gam D(\Gam Ds)+45\Gam^2 D^2s\\
&+225\Gam^3 s
\displaybreak[0]\\
\cD_7s=D^7s&+6D^5(\Gam s)+10D^4(\Gam Ds)+12D^3(\Gam D^2s)\\
&\qquad\qquad+12D^2(\Gam D^3s)+10D(\Gam D^4s)+6\Gam D^5s\\
&+72D^3(\Gam^2 s)+72D^2(\Gam D(\Gam s))+120D^2(\Gam^2 Ds)
+60D(\Gam D^2(\Gam s))\\
&\qquad\qquad+100D(\Gam D(\Gam Ds))+36\Gam D^3(\Gam s)\\
&\qquad\qquad+60\Gam D^2(\Gam Ds)+120D(\Gam^2 D^2s)
+72\Gam D(\Gam D^2s)+72\Gam^2 D^3s\\
&+720D(\Gam^3s)+432\Gam D(\Gam^2s)+432\Gam^2D(\Gam s)+720\Gam^3Ds
\displaybreak[0]\\
\cD_8s=D^8s&+7D^6(\Gam s)+12D^5(\Gam Ds)+15D^4(\Gam D^2s)\\
&\qquad\qquad+16D^3(\Gam D^3s)+15D^2(\Gam D^4s)+12D(\Gam D^5s)+7\Gam D^6s\\
&+105D^4(\Gam^2 s)+112D^3(\Gam D(\Gam s))+192D^3(\Gam^2 Ds)
+105D^2(\Gam D^2(\Gam s))\\
&\qquad\qquad+180D^2(\Gam D(\Gam Ds))+84D(\Gam D^3(\Gam s))
+225D^2(\Gam^2 D^2s)\\
&\qquad\qquad+144D(\Gam D^2(\Gam s))+49\Gam D^4(\Gam s)+84\Gam D^3(\Gam Ds)
+180D(\Gam D(\Gam D^2s))\\
&\qquad\qquad+105\Gam D^2(\Gam D^2s)+192D(\Gam^2 D^3s)+112\Gam D(\Gam D^3s)
+105\Gam^2 D^4s\\
&+1575D^2(\Gam^3s)+1260D(\Gam D(\Gam^2s))+1344D(\Gam^2 D(\Gam s))
+735\Gam D^2(\Gam^2s)\\
&\qquad\qquad+2304D(\Gam^3Ds)+784\Gam D(\Gam D(\Gam s))\\
&\qquad\qquad+735\Gam^2D^2(\Gam s)+1344\Gam D(\Gam^2Ds)+1260\Gam^2D(\Gam Ds)
+1575\Gam^3D^2s\\
&+11025\Gam^4s.
\end{align*}
The combinatorics of the coefficients are simpler than in~\cite{CSS3} and the
numbers are generally smaller; the formulae there are obtained from those here
by expanding the derivatives of $\Gam$ using the product rule. Notice that the
coefficients depend only on the position of the $\Gam$'s, and the coefficients
of the nonlinear terms in $\Gam$ are easily computed as products of the
coefficients of the linear terms, as is clear from the inductive definition of
each $\cD_k$.

We are still free to choose the normalization of the Killing form
$(\cdot\,,\cdot)$: since $-\frac12|\lroot|^2$ is independent of the long root
$\lroot$, we could arrange that this is $1$ and $\Gam=r^D$. This is the
normalization that gives the formulae stated in the introduction for the
conformal case.

\section{Scope of the construction}\label{scope}

We now show that the class of operators constructed in the previous section
includes many standard invariant operators, at least for the `large' parabolic
subgroups occuring in interesting examples. We shall also show that in
conformal geometry, the operators we construct include (at least in the
conformally flat case) those coming from AdS/CFT correspondence for partially
massless fields in string theory.

\subsection{Lagrangian contact structure}
Let us consider the case of a Lagrangian contact structure (see
\cite{C-LC,Slovak,Takeu}).  This is the real split case of the complex
parabolic algebra corresponding to the Dynkin diagram
\begin{equation}
\xbdbx{\alpha_1}{\alpha_2}{}{\alpha_n}{\alpha_{n+1}}\;\;,\; n\geq 1.
\end{equation}
The Lie group is $G=PSL(n+2,\R)$ with Lie algebra $\g=\lie{s}\lie{l}(n+2,\R)$,
$\gr\g$ being equipped with the $|2|$-grading given by block matrices of size
$1,n,1$. The $\g_1$ part decomposes further into a direct sum
$\g_1=\g^L_1\dsum\g^R_1$ of two irreducibles; the $\g_2$ part is
one-dimensional.
 
In geometric terms, we have a contact structure on a real manifold of
dimension $2n+1$ with a direct sum decomposition of the contact distribution
into two Lagrangian subbundles.

To be explicit, we consider the case $n=3$, when the two irreducible
components of $\g_{-1}$ are three dimensional. Let us denote roots
corresponding to both by $e_1,e_2,e_3$ and $f_1,f_2,f_3$ respectively, ordered
so that $e_1,f_1$ are the highest weights for $\g_1$, considered as a
$\p_0$-module, and $e_3,f_3$ are the lowest ones. Let $g$ be the root
corresponding to $\g_{-2}$.
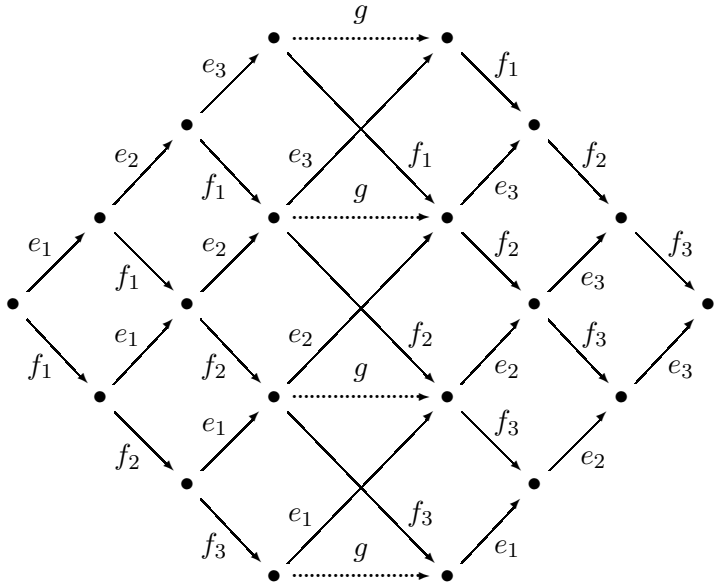
\newarrow{Nul}{}{}{}{}{}
\begin{diagram}[size=1.5em,nohug]
   &      &    &       &    &     &\bullet& &\rDotsto^g
                              & &\bullet&      &    &       &    &    &     \\
   &      &    &       &  &\ruTo^{e_3}&&\rdTo(4,4)&
                         &\ruTo(4,4)&&\rdTo^{f_1}&  &       &    &    &     \\
   &      &    &     &\bullet&      &    &  &
                              &   &    &     &\bullet&      &    &    &     \\
   &      &  &\ruTo^{e_2}&&\rdTo_{f_1}&&\ruNul^{e_3}&        
                       &\rdNul^{f_1}&&\ruTo_{e_3}&&\rdTo^{f_2}&  &    &     \\
   &    &\bullet&      &    &     &\bullet& &\rDotsto^g
                              & &\bullet&      &    &     &\bullet&   &     \\
&\ruTo^{e_1}&&\rdTo_{f_1}&&\ruTo^{e_2}&&\rdTo(4,4)&
                         &\ruTo(4,4)&&\rdTo^{f_2}&&\ruTo_{e_3}&&\rdTo^{f_3} \\
\bullet&  &    &     &\bullet&      &    &  &
                              &   &    &     &\bullet&      &    &  &\bullet\\
&\rdTo_{f_1}&&\ruTo^{e_1}&&\rdTo_{f_2}&&\ruNul^{e_2}&
                       &\rdNul^{f_2}&&\ruTo_{e_2}&&\rdTo^{f_3}&&\ruTo_{e_3} \\
   &    &\bullet&      &    &     &\bullet& &\rDotsto^g
                              & &\bullet&      &    &     &\bullet&   &     \\
   &      &  &\rdTo_{f_2}&&\ruTo^{e_1}&&\rdTo(4,4)&
                         &\ruTo(4,4)&&\rdTo^{f_3}&&\ruTo_{e_2}&  &    &     \\
   &      &    &     &\bullet&      &    &  &
                              &   &    &     &\bullet&      &    &    &     \\
   &      &    &       &  &\rdTo_{f_3}&&\ruNul^{e_1}&
                       &\rdNul^{f_3}&&\ruTo_{e_1}&  &       &    &    &     \\
   &      &    &       &    &     &\bullet& &\rDotsto^g
                              & &\bullet&      &    &       &    &    &
\end{diagram}
The (labelled) Hasse diagram for standard operators is shown above. The labels
on the arrows indicate the `directions' $\lroot$ for the corresponding
operators. We have constructed all operators indicated by full arrows, hence
only the horizontal arrows are missing.

A similar diagram applies in CR geometry, this being another real form of
Lagrangian contact geometry, except that some representations and operators
become conjugate.

\subsection{$G_2$-case}

For a more exotic example, let us consider the split real case of the $G_2$
complex algebra.  In this case, the root system has 12 elements. If we denote
simple roots by $\alpha_1$ (the longer one) and $\alpha_2$ (the shorter one),
then the set of positive roots is
$\{\alpha_1,\alpha_2,\alpha_3,\alpha_4,\alpha_5,\alpha_6\}$ with
$\alpha_3=\alpha_1+\alpha_2,\,\alpha_4=\alpha_1+2\alpha_2,\,
\alpha_5=\alpha_1+3\alpha_2,\, \alpha_6=2\alpha_1+3\alpha_2$.  Let us consider
the case that the parabolic is a Borel (maximal solvable) subgroup of
$G_2$. Then the associated graded algebra $\gr\g$ is $|5|$-graded with
$$
\g_1=\R\cdot\{\al_1\}\dsum\R\cdot\{\al_2\},\,
\g_2=\R\cdot\{\al_3\},\,
\g_3=\R\cdot\{\al_4\},\,
\g_4=\R\cdot\{\al_5\},\,
\g_5=\R\cdot\{\al_6\}.
$$
The (labelled) Hasse diagram then has the form
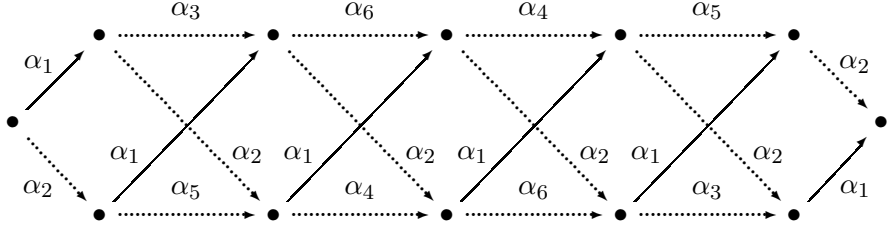
\begin{diagram}[size=1.5em,nohug]
       &             &\bullet&              &\rDotsto^{\al_3}&              
                     &\bullet&              &\rDotsto^{\al_6}&               
                     &\bullet&              &\rDotsto^{\al_4}&              
                     &\bullet&              &\rDotsto^{\al_5}&              
                     &\bullet&              &\\
       &\ruTo^{\al_1}&       &\rdDotsto(4,4)&                &\ruTo(4,4)
                     &       &\rdDotsto(4,4)&                &\ruTo(4,4)
                     &       &\rdDotsto(4,4)&                &\ruTo(4,4)
                     &       &\rdDotsto(4,4)&                &\ruTo(4,4) 
                     &       &\rdDotsto^{\al_2}&\\
\bullet&             &       &              &                &              
                     &       &              &                &              
                     &       &              &                &              
                     &       &              &                &              
                     &       &              &\bullet\\
   &\rdDotsto_{\al_2}&       &\ruNul^{\al_1}&                &\rdNul^{\al_2}
                     &       &\ruNul^{\al_1}&                &\rdNul^{\al_2}
                     &       &\ruNul^{\al_1}&                &\rdNul^{\al_2}
                     &       &\ruNul^{\al_1}&                &\rdNul^{\al_2}
       &       &\ruTo_{\al_1}&\\
       &             &\bullet&              &\rDotsto^{\al_5}&              
                     &\bullet&              &\rDotsto^{\al_4}&              
                     &\bullet&              &\rDotsto^{\al_6}&              
                     &\bullet&              &\rDotsto^{\al_3}&              
                     &\bullet&              &   
\end{diagram}
The operators constructed in section~\ref{ecio}, indicated by full arrows, are
now not so numerous.

There are two other parabolic subgroups of $G_2$ up to isomorphism, one
inducing a $|3|$-grading, the other a $|2|$-grading of the Lie algebra of
$G_2$. For the $|3|$-grading, all roots in $\g_1$ are short, and so no
operators are constructed. For the $|2|$-grading, we have:
$$
\g_1=\R\cdot\{\al_1\}\dsum\R\cdot\{\al_3\}\dsum \R\cdot\{\al_4\}
\dsum\R\cdot\{\al_5\},\,
\g_2=\R\cdot\{\al_6\}.
$$
The Hasse graph in this case is
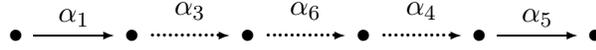
\begin{diagram}[size=2em,nohug]
&\bullet& \rTo^{\al_1} &\bullet& \rDotsto^{\al_3} &\bullet&
\rDotsto^{\al_6} &\bullet& \rDotsto^{\al_4} &\bullet&\rTo^{\al_5} &\bullet
\end{diagram}
with approximately half of the operators constructed in section~\ref{ecio}.

\subsection{The conformal case}\label{conf}
In the even dimensional case, \emph{all} operators in the BGG sequence are
obtained by the construction in section~\ref{ecio}, although there are
nonstandard operators which are not. In the odd dimensional case there is one
arrow in the Hasse diagram which has a special character, as was already noted
in \cite{Gover}. In dimension $2n-1$, $\g=\lie{so}(p,q,\R)$ with $p+q=2n+1$
and a $|1|$-grading. We denote positive roots with root spaces included in
$\g_1$ by $\pm\al_1,\ldots\pm\al_n,\al_{n+1},$ where $\al_{n+1}$ is the short
simple root.  Suppose that the roots $\al_1,\ldots\al_n$ are ordered, i.e.,
that $\al_1$ is the highest among them and $\al_n$ is the smallest.  Then the
(labelled) Hasse diagram has the form
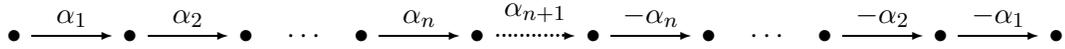
\begin{diagram}[size=2em,nohug]
\bullet& \rTo^{\al_1}& \bullet & \rTo^{\al_2} & \bullet & \cdots &
\bullet& \rTo^{\al_n}& \bullet &\rDotsto^{\al_{n+1}}&\bullet&\rTo^{-\al_n}&
\bullet&\cdots&\bullet& \rTo^{-\al_2}&\bullet&\rTo^{-\al_1}&\bullet
\end{diagram}
The middle operator is labelled by the only short root, hence the construction
of section~\ref{ecio} does not apply; however, the other operators are
all constructed.

\subsection{AdS/CFT for partially massless fields}
The operators constructed in the conformal case include some operators
closely related to the AdS/CFT correspondence for partially massless fields.

We recall that, besides massive or strictly massless fields,
partially massless fields of higher spin were studied on vacuum Einstein
manifolds with a cosmological constant \cite{Buch}. In
Dolan--Nappi--Witten~\cite{DNW}, the following situation was considered.
Let $M$ be an Einstein manifold asymptotic to anti-de Sitter space of
dimension $4$ and let $X$ be its boundary with the induced conformal
structure \cite{GW}.  Then the AdS/CFT correspondence \cite{Malda} yields a
correspondence between a partially massless field $\phi$ on $M$ and a field
$L$ on $X$, satisfying a certain conformally invariant equation (called a
`partial conservation law' in \cite{DNW}).
 
In the case that the field $L$ is a symmetric traceless tensor field
$L^{i_1\ldots i_s}$, the equation is, to leading order,
$$
\nabla_{i_1}\cdots\nabla_{i_{s-n}}L^{i_1\ldots i_s}+\cdots=0.
$$

The source and the target for the equation are easily identified and the
operator is the last operator in the BGG sequence (described in \S\ref{conf}
in the odd dimensional case).  Hence section~\ref{ecio} provides an explicit
formula for a conformally invariant operator of this form.

A detailed study of the case $s=2,n=0$ in \cite{DNW} leads to the equation
on $\R^d$
$$
\nabla_{i}\nabla_{j}L^{ij}+\frac{1}{d-2}R_{ij}L^{ij}=0,
$$
which agrees with the formula of section~\ref{ecio} since the
tracefree parts of $\frac1 {n-2} R_{ij}$ and $\Gam$ agree.

The higher order equations with curvature corrections in section~\ref{ecio}
are natural candidates for the corresponding higher order equations (the cases
with $s-n > 2$) for the field $L$. As an example, let us consider the case
$s=3,n=0$, where we get
\begin{equation}
\nabla_{i}\nabla_{j}\nabla_k L^{ijk}+\frac{2}{d-2}
[\nabla_i(R_{jk}L^{ijk})+R_{jk}\nabla_i L^{ijk}]=0
\end{equation}
for a traceless symmetric tensor field with three indices on $\R^d$.  (Of
course, after expanding the covariant derivatives using the Leibniz rule, we
obtain the formulae already present in~\cite{CSS3,Gover}.)  In the conformally
flat situation, this must be the equation arising from the AdS/CFT
correspondence, because representation theory shows that conformally invariant
operators are unique up to a multiple in this case.  In general, there may be
conformally invariant curvature corrections in the lower order terms.

\appendix
\section{Weyl structures as reductions}

In this appendix we relate our approach to Weyl structures and the original
approach of \v Cap and Slov\'ak~\cite{CSl}, who define a Weyl structure to be
a $P_0$-equivariant section $\sigma$ of $\pi_0\colon\cG\to\cG_0$. For this
definition to make sense, an algebraic Weyl structure $\age$ must be fixed so
that $\pi_0\colon P\to P_0$ is split.

By means of the equation $\pt q(\pt) = \sigma(\pi_0(\pt))$ such a section
$\sigma$ is equivalent to a $P$-invariant trivialization $q\colon
\cG\to\exp\p^\perp$ of the principal $\exp\p^\perp$-bundle
$\pi_0\colon\cG\to\cG_0$, where the $P$-invariance means that $p q(\pt
p)\pi_0(p)^{-1}=q(\pt)$ for all $p\in P,\pt\in \cG$.

On the other hand, using the algebraic Weyl structure $\age$, any geometric
Weyl structure on $M$ is given by $\Ge=(\Ad q)\age$ for a unique $P$-invariant
$q\colon\cG\to \exp \p^\perp$. To summarize:

\begin{prop} Let $(\cG\to M,\cc)$ be a parabolic geometry and fix an algebraic
Weyl structure $\age$. Then there is a natural bijection between Weyl
structures $E$ on $M$ and $P_0$-equivariant sections $\sigma$ of
$\pi_0\colon\cG\to\cG_0$.
\end{prop}

In our development, we used the Weyl structure $\Ge$ to give a $P$-invariant
direct sum decomposition $\sgG\colon \cG\times \T\dsum\p_0\dsum\T\dual\to
\cG\times\g$ and hence write
\begin{equation}\label{decomp}
\sgG^{-1}\circ \cc=\cc_{\T}+\cc_{\p_0}+\cc_{\T\dual}
=\cc_{\T}+\sgG^{-1}\circ \cc_{\p}=\cc_{\T}+\sgG^{-1}\circ \cc_{\Ge}+\rho,
\end{equation}
where $\cc_{\p}$ and $\cc_{\Ge}$ are principal $P$-connections on $\cG$ and
$\rho$ is a $P$-invariant $\p^\perp$-valued horizontal $1$-form inducing the
normalized Ricci curvature $r^D$. On the other hand, an algebraic Weyl
structure $\age$ gives a fixed direct sum decomposition $\age_*\colon
\T\dsum\p_0\dsum\T\dual\to\g$, related to $\sgG$ by conjugating with the
action of $q$, where $\Ge=(\Ad q)\age$.

Since the fixed decomposition is only $P_0$-invariant, \v Cap and Slov\'ak
define the \emph{Weyl form} to be the pull back $\tau=\sigma^*\cc$ of $\cc$ to
$\cG_0$, then decompose $\tau$ into a solder form, a principal
$P_0$-connection and a $P_0$-invariant $\p^\perp$-valued $1$-form. In our
approach~\eqref{decomp} is a $P$-invariant lift of this decomposition to
$\cG$.

\section{Dependence of $D$ and $r^D$ on the Weyl structure}

In equation~\eqref{D1-change}, we obtained the (infinitesimal) dependence of
the Ricci-corrected Weyl connection $D^{(1)}$ on the Weyl structure. We now do
the same for the Weyl connection $D$ and the normalized Ricci curvature $r^D$
(we do not need these results in the body of the paper, but include them for
general interest).

\begin{prop} For $\ag\in \Cinf(M,\pM^\perp)$ and $X\in TM$, $\pd_\ag r^D(X)
= -D\low_X\ag +[\ag,X]^{\gem}_{\smash{\pM^\perp}}$, where $X$ is lifted to
$\gM$ and the Lie bracket is projected onto $\pM^\perp$ using $\gem$.
\proofof{prop} $r^D$ is the $\pM^\perp$-valued $1$-form on $M$ induced by
$\rho=\Ge^*\MCw+\cc_{\T\dual}$, where $\Ge^*\MCw_y=\MCw_{\Ge(y)}d\Ge_y$ and
$\cc_{\T\dual}$ is shorthand for the $\T\dual$ component
$(\sgG^{-1}\cc)_{\T\dual}$.  Viewing $\ag$ as a $P$-invariant
$\p^\perp$-valued function on $\cG$ with $\pd_\ag\Ge=-\ag$, we easily compute
that $\pd_\ag(\Ge^*\MCw)=-d\ag-[\ag,\Ge^*\MCw]$ and $\pd_\ag\cc_{\T\dual}=
(\sgG^{-1}[\ag,\cc])\low_{\T\dual}- [\ag,\cc_{\T\dual}]=
[\ag,\sgG\cc_{\p_0}]+(\sgG^{-1}[\ag,\sgG \cc_{\T}])\low_{\T\dual}$. Hence
\begin{equation*}
\pd_\ag\rho = -(d\ag-[\Ge^*\MCw,\ag] +[\sgG \cc_{\p_0},\ag])
 +[\ag,\cc_{\T}]^{\Ge}_{\p^\perp}
\end{equation*}
as required, since $\cc_{\Ge}=-\Ge^*\MCw+\sgG\cc_{\p_0}$ is the principal
connection inducing $D$.
\end{proof}

\begin{prop} Let $\vE$ be a filtered $P$-bundle and $\vs$ a section of $\vE$.
Then for $\ag\in \Cinf(M,\pM^\perp)$ and $X\in TM$, $\pd_\ag D\low_X\vs=
([\ag,X]^{\gem}_{\smash{\p\low_{M,0}}}+D\low_X\ag)\act\vs$.
\proofof{prop} $D \vs = \gem^{-1}_\vE D^{(1)} (\gem_\vE \vs)$ and so
$\pd_\ag D\low_X \vs =\gem^{-1}_\vE ([\ag,X]^{\gem}_{\smash{\pM}}\act
(\gem_\vE \vs)) +D\low_X(\ag\act\vs)-\ag\act D\low_X\vs$. As $\gr\vE$ is
semisimple, the result follows.
\end{proof}
Since $D^{(1)}_X\vs=D\low_X \vs+r^D(X)\act\vs$, these computations are not
independent: we check
\begin{equation*}
\pd_\ag(D\low_X \vs+r^D(X)\act\vs)
= [\ag,X]^{\gem}_{\smash{\p\low_{M,0}}}\act\vs
+[\ag,X]^{\gem}_{\smash{\pM^\perp}}\act\vs
= [\ag,X]^{\gem}_{\smash{\pM}}\act\vs
\end{equation*}
in accordance with equation~\eqref{D1-change}.  Unlike $D^{(1)}$, the Weyl
connection $D$ does not depend algebraically on the Weyl structure.  This, of
course, was the whole reason for introducing Ricci corrections in the first
place.

%
\newcommand{\bauth}[1]{\mbox{#1}, \ignorespaces}
\newcommand{\bart}[1]{\textit{#1}, \ignorespaces}
\newcommand{\bjourn}[3]{#1 \textbf{#2} (#3) \ignorespaces}
\newcommand{\bbook}[1]{\textsl{#1}, \ignorespaces}
\newcommand{\bseries}[2]{#1 \textbf{#2}, \ignorespaces}
\newcommand{\bpp}[2]{#1--#2.}

\end{document}